\documentclass[review]{elsarticle}
\usepackage{amsopn}

\newcommand{\IR}{\mathbb{R}}

\newcommand{\CC}{\mathcal{C}}

\newcommand{\IE}{\mathscr{I}}
\newcommand{\MM}{\mathcal{M}}
\newcommand{\IM}{\mathscr{M}}
\newcommand{\IV}{\mathscr{V}}

\newcommand{\IF}{\mathscr{F}}

\newcommand{\HH}{\eta}

\newcommand{\Domain}{\mathscr{D}}
\newcommand{\EQ}[1]{(\ref{eq:#1})}
\newcommand{\SEC}[1]{\ref{sec:#1}}
\newcommand{\FIG}[1]{\ref{fig:#1}}
\newcommand{\feq}{f_{\mathrm{eq}}}
\newcommand{\llbrack}{{\lbrack\hspace{-1.25pt}\lbrack}}
\newcommand{\rrbrack}{{\rbrack\hspace{-1.25pt}\rbrack}}

\newcommand{\jump}[1]{\llbrack #1 \rrbrack}

\newcommand{\Kn}{\mathrm{Kn}}

\newcommand*{\defeq}{\mathrel{\rlap{%
                     \raisebox{0.3ex}{$\m@th\cdot$}}%
                     \raisebox{-0.3ex}{$\m@th\cdot$}}%
                     =}
\newcommand*{\eqdef}{=\mathrel{\rlap{%
                     \raisebox{0.3ex}{$\m@th\cdot$}}%
                     \raisebox{-0.3ex}{$\m@th\cdot$}}%
                     }

\DeclareMathOperator{\argmin}{arg\:min}

\usepackage{lmodern}
\usepackage{amsmath}
\usepackage{amsthm}
\usepackage{mathrsfs}
\usepackage{amssymb}
\usepackage{bm}
\usepackage{subfig}
\usepackage[margin=2.5cm]{geometry}
\usepackage{color}
\usepackage{epstopdf}
\usepackage{citesort}
\usepackage{algpseudocode}
\usepackage{boxedminipage}
\usepackage{float}
\usepackage{url}
\newfloat{Algorithm}{t}{alg}

\newtheorem*{remark}{Remark}
\begin{document}

\begin{frontmatter}

\title{Error estimation and adaptive moment hierarchies for goal-oriented approximations of the Boltzmann equation} 

\author{M.R.A. Abdelmalik\corref{mycorrespondingauthor}}  %
\ead{m.abdel.malik@tue.nl}
\author{E.H. van Brummelen}
\address{Department of Mechanical Engineering, Eindhoven University of Technology, Netherlands}
\cortext[mycorrespondingauthor]{Corresponding author}
\begin{abstract}
This paper presents an a-posteriori goal-oriented error analysis for a numerical approximation of the steady Boltzmann equation based on a moment-system approximation in velocity dependence and a discontinuous Galerkin finite-element (DGFE) approximation in position dependence. We derive computable error estimates and bounds for general target functionals of solutions of the steady Boltzmann equation based on the DGFE moment approximation. The a-posteriori error estimates and bounds are used to guide a model adaptive algorithm for optimal approximations of the goal functional in question. We present results for one-dimensional heat transfer and shock structure problems where the moment model order is refined locally in space for optimal approximation of the heat flux.
\end{abstract}


\end{frontmatter}

\section{Introduction}
\label{sec:intro}
The Boltzmann equation provides a description of the molecular dynamics of fluid flows based on their one-particle phase-space distribution. The Boltzmann equation encapsulates all conventional macroscopic flow models in the sense that its limit solutions correspond to solutions of the compressible Euler and Navier-Stokes equations \cite{Bardos:1991vf,Esposito:1994ca}, the incompressible Euler and Navier--Stokes equations~\cite{Golse:2004oe,Lions:2001sj}, the incompressible Stokes equations \cite{Lions:2001wb} and the incompressible Navier-Stokes-Fourier system \cite{Levermore:2010kx}; see \cite{Saint-Raymond2009} for an overview. Fluid flow problems generally exhibit locally varying deviations from a continuum description. Therefore, the Boltzmann equation is uniquely suited to describe flows with varying rarefaction regimes. Applications in which rarefaction effects play a significant role are multitudinous, including gas flow problems involving large mean free paths in high-altitude flows and hypobaric applications such as chemical vapor deposition; see \cite{Cercignani2000,Struchtrup2005} and references therein for further examples. Moreover, the perpetual trend toward miniaturization in science and technology renders accurate descriptions of fluid flows in the transitional molecular/continuum regime of fundamental technological relevance, for instance, in nanoscale applications, micro-channel flows or flow in porous media~\cite{shen2005}. The Boltzmann equation also provides a prototype for kinetic models in many other applications that require a description of the collective behavior of large ensembles of small particles, for instance, in semi-conductors \cite{jungel2009}, in plasmas and fusion and fission devices \cite{miyamoto2004}  and in dispersed-particle flows such as in fluidized-bed reactors \cite{reeks1991,reeks1992,reeks1993}.

Numerical approximation of the Boltzmann equation poses a formidable challenge, on account of its high dimensional setting: for a problem in~$D$ spatial dimensions, the one-particle phase-space is~$2D$ dimensional. The corresponding computational complexity of conventional discretization methods for (integro\nobreakdash-)differential equations, such as finite-element methods with uniform meshes, is prohibitive. Numerical approximations of the Boltzmann equation have been predominantly based on particle methods, such as the Direct Simulation Monte Carlo (DSMC) method~\cite{Bird:1970,Bird:1994}. Convergence proofs for these methods~\cite{Wagner1992} however convey that their computational complexity depends sensitively on the Knudsen number, and the computational cost becomes prohibitive in the fluid-dynamical limit. Moreover, from an approximation perspective, DSMC can be inefficient, because it is inherent to the underlying Monte-Carlo process that the approximation error decays only as $n^{-\frac{1}{2}}$ as the number of simulation molecules, $n$, increases; see, for instance, \cite[Thm.~5.14]{Klenke2008}. Efficient computational modeling of fluid flows in the transitional molecular/continuum regime therefore remains an outstanding challenge.

An alternative approximation technique for the Boltzmann equation is the method of moments \cite{Grad1949,Levermore1996,Struchtrup2005,Abdelmalik2016a}. The method of moments represents a general statistical approximation technique which identifies parameters of an approximate distribution based on its moments \cite{Matyas1999}. Application of the method of moments to the Boltzmann equation engenders a system of evolution equations for the moments (weighted averages) of the phase-space distribution. In \cite{Abdelmalik2016a,Abdelmalik2016b} it was shown that moment systems can be conceived of as Galerkin approximations, in velocity dependence, of the Boltzmann equation in renormalized form, where larger moment systems correspond to refined subspaces. 
By virtue of the hierarchical structure of the considered subspaces, moment systems form a hierarchy of models that bridge the transitional molecular/continuum flow regime. For suitable choices of the renormalization map,
moment systems are symmetric hyperbolic and are well-posed in the corresponding sense~\cite{Majda:1984wj}.
The symmetric hyperbolic structure of moment systems moreover allows one to devise stable Galerkin approximations 
in position dependence \cite{Barth2006,hartmann2002}.

By virtue of their Galerkin form and their inherent hierarchical structure, moment systems are 
ideally suited for (goal-oriented)
model adaptivity. Goal-oriented a-posteriori error-estimation and adaptivity exploit the Galerkin method to construct finite-element spaces that yield an optimal approximation of a particular functional of the solution (goal or target functional). Goal-oriented methods for mesh refinement date back to the pioneering work of Becker and Rannacher \cite{becker1996,becker2001}, Oden and Prudhomme \cite{oden2001}, Giles and S{\"u}li \cite{giles2002}, Houston and Hartmann \cite{hartmann2002,hartmann2003}, and Houston and S{\"u}li \cite{houston2001}. In addition to mesh adaptivity, goal-oriented approaches have also been considered in the context of model adaptivity where the algorithm adapts between a coarse and a fine model, such that the fine model is only applied in regions of the domain that contribute most significantly to the goal functional in question. For examples of goal-oriented model adaptivity, we refer to \cite{oden2000} for application of goal-oriented model adaptivity to heterogeneous materials, to \cite{bauman2009} for goal-oriented atomistic/continuum adaptivity in solid materials, and to \cite{opstal2015} for goal-oriented adaptivity between the Stokes equations and the Navier-Stokes equations. The Galerkin form of moment methods enables the construction of accurate a-posteriori error estimates, while
the hierarchical structure provides an intrinsic mode of refinement. In principle, the error estimate can serve to
guide simultaneous (anisotropic) mesh and moment refinement. If the adaptive strategy is restricted to moment refinement
only, i.e. the finite-element approximation in position dependence is fixed, the adaptive procedure can be 
viewed as a goal-oriented model-adaptive strategy~\cite{Ainsworth2011} that adapts between the models
in the moment-system hierarchy to construct an optimal approximation to the goal function in question.

The purpose of this paper is to derive a-posteriori error estimates, measured in terms of a certain target functional, for a position-velocity Galerkin approximation of the steady Boltzmann equation. The Galerkin approximation is based on a moment-system approximation in velocity dependence and a discontinuous Galerkin finite element (DGFE) approximation in position dependence. We propose a goal-adaptive refinement procedure that locally adapts the order of the moment system to locally resolve rarefaction effects corresponding to their contribution to the error in the quantity of interest. By only targeting regions with the largest contributions to the error the model adaptive strategy yields an approximation that is quasi-optimal for the goal-functional.

The proposed adaptive moment method can alternatively be classified as a heterogeneous multiscale method of type~A; 
see~\cite{e2003}. Multiscale methods of type~A introduce a decomposition of the spatial domain into a subset 
where a {\em macroscopic\/} (or {\em coarse\/}, {\em simple\/}) model suffices, and a complementary subset where 
a {\em microscopic\/} (or {\em fine\/}, {\em sophisticated\/}) description is required. The proposed adaptive moment method introduces an element-wise domain decomposition strategy where, locally, different levels of the moment hierarchy are used to approximate the solution to the Boltzmann equation. The goal-adaptive algorithm provides an automated strategy for model refinement such that an optimal approximation of the solution of the Boltzmann equation is obtained for the goal 
functional under consideration.

The remainder of this paper is organized as follows. Section~\ref{sec:BoltzMoms} surveys standard structural properties of the Boltzmann equation that are retained in the moment-system approximation. Section \ref{sec:Galerkin} introduces the Galerkin approximation of the stationary Boltzmann equation. We will present the moment system as a Galerkin approximation, in velocity dependence, of the Boltzmann equation in renormalized form. In space dependence we discretize the Boltzmann equation using the discontinuous Galerkin finite-element method. Section \ref{sec:ErrEst} presents the derivation of a-posteriori error estimates for the DGFE moment approximation. In Section~\ref{sec:AdaptAlg} we devise an adaptive algorithm for the steady DGFE moment method that exploits inter-element cancellation errors. In Section \ref{sec:NumRes} we apply the adaptive algorithm to the heat transfer and shock structure Riemann problems \cite{Courant1999,Toro2013}. Finally, Section~\ref{sec:Conc} presents a concluding discussion.

\section{The Boltzmann equation}
\label{sec:BoltzMoms}
Consider a monatomic gas contained in a fixed spatial domain $\Omega\subset\mathbb{R}^D $. Kinetic theory describes the state of such a gas by a non-negative (phase-space) density $f=f(t,\bm x,\bm v)$ over the single-particle phase space $\Omega\times\mathbb{R}^D$. The evolution of $f$ is considered to be governed by the Boltzmann equation,
\begin{align}\label{eq:Boltzmann}
\partial_t f + v_i\partial_{x_i} f=\CC(f)
\end{align}
where the summation convention applies to repeated indices and the collision operator $f\mapsto\CC(f)$ acts only on the $\bm v=(v_1,\ldots,v_D)$ dependence  of $f$ locally at each $(t,\bm x)$. The collision operator is assumed to possess certain conservation, symmetry and dissipation properties, viz., conservation of mass, momentum and energy, Galilean invariance and dissipation of appropriate entropy functionals. Moreover, it is assumed that the collision operator exhibits certain positivity properties. These fundamental properties of the collision operator have been treated in detail in \cite{Abdelmalik2016a,Abdelmalik2016b} and are merely repeated here for completeness and coherence.

To elaborate the conservation properties of the collision operator, let $\langle \cdot \rangle$ denote integration in the velocity dependence of any scalar, vector or matrix valued measurable function over $D$\nobreakdash-dimensional Lebesgue measure. A function $\alpha:\IR^D\to\IR$ is called a {\em collision invariant\/} of $\CC$ if
\begin{align}
\label{eq:DefCollInvariant}
\langle \alpha\,\CC(f) \rangle = 0  \qquad \forall f\in \Domain(\CC),
\end{align}
where $\Domain(\CC)\subset{}L^1(\IR^D,\IR_{\geq{}0})$ denotes the domain of~$\CC$, which we consider to be a subset of the almost everywhere nonnegative
Lebesgue integrable functions on~$\IR^D$.
Equation~\EQ{DefCollInvariant} associates a scalar conservation law with each collision invariant:
\begin{equation}
\label{eq:ConsvLaw}
\partial_t\langle \alpha f\rangle+\partial_{x_i} \langle v_i \alpha f \rangle = 0
\end{equation}
We insist that $\{1,v_1,\ldots,v_D,|\bm v|^2\}$ are collision invariants of $\CC$ and that
the span of this set contains all collision invariants, i.e.
\begin{equation}\label{eq:CollInvariant}
\langle{}\alpha\,\CC(f)\rangle=0\quad\forall{}f\in\Domain(\CC)
\quad\Leftrightarrow\quad
\alpha\in\mathrm{span}\{1,v_1,\ldots,v_D,|\bm v|^2\}=:\IE.
\end{equation}
The moments $\langle{}f\rangle$,
$\langle{}v_if\rangle$ and $\langle{}|{\bm v}|^2f\rangle$, correspond to the mass-density, the (components of) momentum-density and the energy-density, respectively.
Accordingly, the conservation law~\EQ{ConsvLaw} implies that~\EQ{Boltzmann} conserves mass, momentum and energy.

The assumed symmetry properties of the collision operator pertain to commutation with translational and rotational transformations.
In particular, for all vectors $\bm u\in\mathbb{R}^D$ and all orthogonal tensors $\mathcal{O}:\IR^D\to\IR^D$, we define the translation transformation
$T_{\bm u}:\Domain(\CC)\to\Domain(\CC)$ and the rotation transformation $T_{\mathcal{O}}:\Domain(\CC)\to\Domain(\CC)$ by:
\begin{alignat}{2} \label{eq:Tran}
(T_{\bm u}f)(\bm v)&=f(\bm u-\bm v)&\qquad&\forall{}f\in\Domain(\CC)
\\
\label{eq:Rot}
(T_{\mathcal{O}}f)(\bm v)&=f(\mathcal{O}^*\bm v) &\qquad&\forall{}f\in\Domain(\CC)
\end{alignat}
with $\mathcal{O}^*$ the Euclidean adjoint of $\mathcal{O}$. Note that the above transformations act on the $\bm v$-dependence only.
It is assumed that $\CC$ possesses the following symmetries:
\begin{equation}
\label{eq:GalilInvar}
\CC(T_{\bm u}f)=T_{\bm u}\CC(f),\qquad
\CC(T_{\mathcal{O}}f)=T_{\mathcal{O}}\CC(f)
\end{equation}
The symmetries (\ref{eq:GalilInvar}) imply that (\ref{eq:Boltzmann}) complies with Galilean invariance, i.e.
if $f(t,\bm x,\bm v)$ satisfies the Boltzmann equation~\EQ{Boltzmann}, then for arbitrary ${\bm u}\in\IR^D$ and arbitrary orthogonal $\mathcal{O}:\IR^D\to\IR^D$, so
do $f(t,\bm x - {\bm u}t,{\bm v}-{\bm u})$ and $f(t,\mathcal{O}^*{\bm x},\mathcal{O}^*{\bm v})$.

The entropy dissipation property of~$\CC$ is considered in the extended setting of~\cite[Sec.~7]{Levermore1996}, from which we
derive the following definition: a convex function $\HH:\IR_{\geq{}0}\to\IR$ is called an {\em entropy for $\CC$\/} if
\begin{equation}
\label{eq:Dissipation}
\langle \HH'(f)\,\CC(f) \rangle \leq 0, \qquad \forall f\in\Domain(\CC)
\end{equation}
with $\HH'(f)$ the derivative of $\HH(f)$, and if for every $f\in\Domain(\CC)$ the following equivalences hold:
\begin{equation}
\label{eq:Equilibrium}
\CC(f) = 0
\quad\Leftrightarrow\quad
\langle \HH'(f)\,\CC(f)\rangle  = 0
\quad\Leftrightarrow\quad
\HH'(f)\in\IE
\end{equation}
Relation (\ref{eq:Dissipation}) implies that~$\CC$ dissipates the local entropy density $\langle\HH(\cdot)\rangle$, which leads to an
abstraction of Boltzmann's H-theorem for~\EQ{Boltzmann}, asserting that solutions of the Boltzmann equation (\ref{eq:Boltzmann}) satisfy the
local entropy-dissipation law:
\begin{align}
\label{eq:EntDiss}
 \partial_t\langle \HH(f) \rangle+\partial_{x_i} \langle v_i \HH(f) \rangle = \langle \CC(f)\, \HH'(f)\rangle \leq 0\,.
\end{align}
The functions $\langle \HH(f) \rangle$,  $\langle v_i \HH(f) \rangle$ and $\langle \HH'(f)\, \CC(f)\rangle$ are referred to as entropy density, entropy flux and entropy-dissipation rate, respectively.
The first equivalence in~(\ref{eq:Equilibrium}) characterizes local equilibria of~$\CC$
by vanishing entropy dissipation, while the second equivalence indicates the form of such local equilibria.
For spatially homogeneous initial data, $f_0$, equations~\EQ{Dissipation} and~\EQ{Equilibrium} suggest that equilibrium solutions, $\feq$,
of~\EQ{Boltzmann} are determined by:
\begin{equation}
\label{eq:LegEq}
\feq=\argmin\big\{\langle\HH(f)\rangle:f\in\Domain(\CC),\langle\alpha f \rangle=\langle{}\alpha f_0\rangle\:\:\forall\alpha\in\IE\},
\end{equation}
Equation~\EQ{LegEq} identifies equilibria as minimizers\footnote{We adopt the sign convention of diminishing entropy.} of the entropy, subject to the constraint that the invariant moments are identical to the invariant moments of the initial distribution.

The standard definition of entropy corresponds to a density $f\mapsto{}\langle{}f\log{}f+f\alpha\rangle{}$ where $\alpha\in\IE$ is any collision invariant. The corresponding local equilibria of $\CC(f)$ defined by (\ref{eq:Equilibrium}) are characterized by Maxwellians $\MM$, i.e. distributions of the form
\begin{equation} \label{eq:Maxwellian}
\MM({\bm v}):=
\MM_{(\varrho,{\bm u},T)}({\bm v}) :=
\frac{\varrho}{(2\pi{}RT)^{\frac{D}{2}}}\exp\left(-\frac{|\bm v-\bm u|^2}{2RT}\right)
\end{equation}
for some $(\varrho, {\bm u}, T)\in\mathbb{R}_{>0}\times\mathbb{R}^D\times\mathbb{R}_{>0}$ and a certain gas constant $R\in\IR_{>0}$. 

We admit distributions that vanish on sets with nonzero measure. To accommodate such distributions, we introduce an auxiliary non-negativity condition on the collision operator, in addition to~(\ref{eq:Dissipation}) and~(\ref{eq:Equilibrium}). The non-negativity condition insists that $\CC(f)$ cannot be negative on zero sets of~$f$:
\begin{equation}
\label{eq:NonNeg}
\CC(f)\big|_{\mathrm{supp}^c(f)}\geq{}0
\end{equation}
where $\mathrm{supp}^c(f)$ denotes the zero set of~$f$, i.e. the complement in $\IR^D$ of the closed support of~$f$. Condition~(\ref{eq:NonNeg}) encodes that the collision operator cannot create locally negative distributions. It can be verified that~(\ref{eq:NonNeg}) holds for a wide range of collision operators,
including the BGK operator~\cite{Bhatnagar:1954hc}, the multi-scale generalization of the BGK operator introduced in~\cite{Levermore1996}, and all collision operators that are characterized by a (non-negative) collision kernel.

\section{Galerkin approximation of the Boltzmann equation}
\label{sec:Galerkin}
In this section we derive a position-velocity Galerkin approximation for the steady Boltzmann equation in renormalized form. The Galerkin approximation is based on a moment-system approximation in velocity dependence and a discontinuous Galerkin approximation in position dependence. The moment system approximation and its equivalence to a Galerkin approximation have been presented in \cite{Abdelmalik2016a,Abdelmalik2016b} and are repeated here for completeness. For the space DGFE moment approximation we will use the numerical flux derived in~\cite{Abdelmalik2016b}. Although we restrict ourselves to steady problems in this work, for transparency, we present the moment formulation in the time-dependent setting.

\subsection{Velocity discretization using moment-system hierarchies}
Our semi-discretization of the Boltzmann equation with respect to the velocity dependence is based on velocity moments of the one-particle marginal. These velocity moments are defined over $\mathbb{R}^D$, and therefore we regard finite dimensional approximations of $f(t,\bm x, \bm v)$ in (\ref{eq:Boltzmann}) that are integrable over $\mathbb{R}^D$ in velocity dependence. To that end, we consider a Galerkin subspace approximation of the Boltzmann equation in renormalized form, where the renormalization maps to integrable functions. To elucidate the renormalization, let $\IM$ denote an $M$\nobreakdash-dimensional subspace of $D$-variate polynomials and let $\{m_i(\bm v)\}_{i=1}^M$  represent a corresponding basis. We consider the renormalization map $\beta:\IM\to\IF$, where 
\begin{multline}
\IF:=\big\{f\in{}\Domain(\CC): f\geq0,
 mf\in{L}^1(\mathbb{R}^{D}),\\
 \bm v mf\in{L}^1(\mathbb{R}^{D},\IR^D),\,
 m\CC(f)\in{}L^1(\IR^D)\:
 \forall m\in\IM\big\}.
\end{multline}
The moment system can then be written in the Galerkin form:
\begin{multline}  \label{eq:GalLevMomCls}
  \text{\it Find }
  g\in{}\mathscr{L}\big((0,T)\times\Omega;\IM\big):\\
  \partial_t\big\langle m \beta(g)\big\rangle+\partial_{x_i}\big\langle
mv_i\beta(g)\big\rangle = \big\langle m\CC(\beta{}(g))\big\rangle
  \\
  \forall m\in{}\IM\text{ a.e. }(t,\bm x)\in(0,T)\times\Omega
\end{multline}
where $\mathscr{L}\big((0,T)\times\Omega;\IM\big)$ represents a suitable vector space of functions from the considered time interval $(0,T)$ and spatial domain~$\Omega$ into~$\IM$.
The symmetry and conservation properties are generally retained in (\ref{eq:GalLevMomCls}) by a suitable selection of the subspace~$\IM$, namely that $\IM$ contains the collision invariants $\IE$, and is closed under the actions of $T_{\bm u}$ and $T_{\mathcal{O}}$; cf.~(\ref{eq:Tran}) and~\EQ{Rot}. To retain entropy dissipation as in (\ref{eq:EntDiss}), we consider a renormalization map and entropy function pair $\{\beta,\eta\}$ that are related by $\beta^{-1}(\cdot)=\eta'(\cdot)$. Entropy dissipation then follows directly from Galerkin orthogonality in (\ref{eq:GalLevMomCls}):
\begin{align}
  \partial_t\big\langle \eta(\beta(g))\big\rangle+\partial_{x_i}\big\langle
  v_i\eta(\beta(g))\big\rangle = \big\langle \eta'(\beta(g))\CC(\beta{}(g))\big\rangle\leq0;
\end{align}
see~\cite{Abdelmalik2016a,Abdelmalik2016b} for more details.

We consider a family of renormalization maps and corresponding entropy functions
according to
\begin{equation}
\label{eq:phimap}
\beta(g)=\mathcal{B}\left(1+\frac{g}{N}\right)^N_+ 
\qquad 
\eta(f) = f\bigg(\frac{N^2}{1+N}\bigg(\frac{f}{\mathcal{B}}\bigg)^{1/N}-N\bigg) + \mathcal{B}\frac{N}{1+N}
\end{equation}
where $(\cdot)_+:=\frac{1}{2}(\cdot)+\frac{1}{2}|\cdot|$, $N$ is a positive integer and $\mathcal{B}$ is some suitable background distribution; see also~\cite{Abdelmalik2016a}. One can infer that indeed $\beta^{-1}(\cdot)=\eta'(\cdot)$
and that $\eta$ is strictly convex on~$\IR_{>0}$. The entropy function in~\EQ{phimap} corresponds to a relative
entropy associated with a $\varphi$\nobreakdash-divergence~\cite{Csiszar1972} with respect to the background 
measure~$\mathcal{B}$. In particular, it holds that $\eta(f)=\mathcal{B}\,\varphi(f/\mathcal{B})$ with $\varphi$ according to:
\begin{equation}
\label{eq:phidivergence}
\varphi(\cdot)=(\cdot)\bigg(\frac{N^2}{1+N}(\cdot)^{1/N}-N\bigg) + \frac{N}{1+N}
\end{equation}
The renormalization map $g\mapsto\beta(g)$ corresponds to a divergence-based moment-closure relation in 
the sense that it associates the following distribution with a given moment vector~$\bm{\mu}\in\IR^M$:
\begin{equation}
\argmin\big\{\langle\mathcal{B}\,\varphi(f/\mathcal{B})\rangle: \langle \bm{m} f \rangle=\bm\mu\big\}
\end{equation}
i.e. the closure relation minimizes the divergence-based relative entropy subject to the constraint that
its moments coincide with the given moments~$\bm{\mu}$.
It was shown in~\cite{Abdelmalik2016a} that the renormalization map~(\ref{eq:phimap})
engenders well-posed moment systems.

\begin{remark}
\label{rem:remark1}
Adoption of a $\varphi$-divergence-based entropy stipulates that this entropy
satisfies~\EQ{Dissipation} and~\EQ{Equilibrium} for a meaningful class of collision
operators subject to~(\ref{eq:NonNeg}). In~\cite{Abdelmalik2016a} it is has been
shown that the class of admissible collision operators includes the BGK
operator~\cite{Bhatnagar:1954hc} and the multi-scale generalization of the BGK
operator in~\cite{Levermore1996}.
\end{remark}

\begin{remark}
It is noteworthy that in the limit $N\rightarrow\infty$, the renormalization map and corresponding entropy 
in~(\ref{eq:phimap}) recover Levermore's moment-closure relation~\cite{Levermore1996}, viz.~$\mathcal{B}\exp(g)$, and the 
Kullback-Leibler divergence~\cite{Kullback1951} relative to $\mathcal{B}$, viz. $\langle f\log(f/\mathcal{B})\rangle$,
respectively; see \cite{Abdelmalik2016a} for more details.

\end{remark}

\subsection{The DGFE moment approximation}
\label{sec:DGFEM}
For the semi-discretization of (\ref{eq:GalLevMomCls}) with respect to the position dependence, we consider the discontinuous Galerkin finite-element method~\cite{Di-Pietro:2012kx}. Henceforth we restrict ourselves to the stationary problem corresponding to (\ref{eq:GalLevMomCls}). Let $\mathcal{H}:=\{h_1,h_2,\ldots\}\subset\IR_{>0}$ denote a strictly decreasing sequence of mesh parameters whose only accumulation point is~$0$. Consider a corresponding mesh sequence $\mathcal{T}^{\mathcal{H}}$, viz., a sequence of covers of the domain by non-overlapping element domains $\kappa\subset\Omega$. We impose on $\mathcal{T}^{\mathcal{H}}$ the standard conditions of regularity, shape-regularity and quasi-uniformity with respect to~$\mathcal{H}$; see, for instance, \cite{Di-Pietro:2012kx} for further details. To introduce the DGFE approximation space, let $\mathcal{P}_p(\kappa)$ denote the set of $D$-variate polynomials of degree at most $p$ in an element domain $\kappa\subset\mathbb{R}^D$. For any $h\in\mathcal{H}$, we indicate by $V^{h,p}(\Omega)$ the DGFE approximation space:
\begin{equation}\label{eq:DGspace}
V^{h,p}(\Omega)=\{g\in{}L^2(\Omega):\ g|_{\kappa}(t,\bm{x})\in\mathcal{P}_p(\kappa), \ \forall\kappa\in\mathcal{T}^h\},
\end{equation}
and by $V^{h,p}(\Omega,\IM)$ the extension of $V^{h,p}(\Omega)$ to $\IM$-valued functions:
\begin{equation}
V^{h,p}(\Omega,\IM)=V^{h,p}(\Omega)\times\IM=\{\lambda_1m_1+\cdots+\lambda_Mm_M:\lambda_i\in{}V^{h,p}(\Omega)\}
\end{equation}
Let us note that for simplicity we have assumed that the dimension of the moment space $M=\dim(\mathscr{M})$ 
is uniform on~$\mathcal{T}^h$. However, this assumption is non-essential and can be dismissed straightforwardly, 
i.e. the (dimension of the) moment space can be selected element-wise. This in fact provides the basis for the 
model-adaptive strategy 
in~Section \ref{sec:AdaptAlg}, which assigns different moment orders $M_{\kappa}$ to the elements $\kappa\in\mathcal{T}^h$.

To facilitate the presentation of the DGFE formulation, we introduce some further notational conventions. For any $h\in\mathcal{H}$, we indicate by $\mathcal{I}^h=\{\mathrm{int}(\partial\kappa\cap\partial\hat{\kappa}):\kappa,\hat{\kappa}\in\mathcal{T}^h,\kappa\neq\hat{\kappa}\}$ the collection of inter-element edges, by $\mathcal{B}^h=\{\mathrm{int}(\partial\kappa\cap\partial\Omega):\kappa\in\mathcal{T}^h\}$ the collection of boundary edges and by $\mathcal{S}^h=\mathcal{B}^h\cup\mathcal{I}^h$ their union. With every edge we associate a unit normal vector~$\bm{\nu}^e$. The orientation of~$\bm{\nu}^e$ is arbitrary except on boundary edges where $\bm{\nu}^e=\bm{n}|_e$. For all interior edges, let $\kappa_{\pm}^e\in\mathcal{T}^h$ be the two elements adjacent to the edge~$e$ such that the orientation of $\bm{\nu}_{\pm}^e=\pm\bm{\nu}^e$ is exterior to~$\kappa_{\pm}^e$. We define the edge-wise jump operator according~to:
\begin{equation}
\label{eq:jumpmean}
\jump{v}=
\begin{cases}
(v_+\bm{\nu}_++v_-\bm{\nu}_-)&\text{ if }e\in\mathcal{I}^h
\\
(v-v_B)\bm{\nu}^e&\text{ if }e\in\mathcal{B}^h
\end{cases}
\end{equation}
where $v_+$ and $v_-$ refer to the restriction of the traces of $v|_{\kappa_+}$ and $v|_{\kappa_-}$ to~$e$. To derive the  DGFE formulation of the closed moment system (\ref{eq:GalLevMomCls}), we note that for any $\psi\in V^{h,p}(\Omega,\IM)$ there holds
\begin{equation}
\label{eq:DGweight}
\sum_{\kappa\in\mathcal{T}^h}\int_{\kappa}\langle
\psi\,\partial_{x_i} v_i
\beta(g)\rangle
=
\sum_{\kappa\in\mathcal{T}^h}\int_{\kappa}
\langle \psi\,\CC(\beta(g))\rangle
\end{equation}
Using the product rule and integration by parts, \EQ{DGweight} can be reformulated in weak form. The left member of~\EQ{DGweight} can be recast into
\begin{equation}
\begin{split}
\label{eq:DGform2}
\sum_{\kappa\in\mathcal{T}^h}\int_{\kappa}\langle
\psi\, \partial_{x_i} v_i
\beta(g)\rangle
&=
\sum_{\kappa\in\mathcal{T}^h}
\int_{\partial\kappa}\langle
\psi \,
 v_i\, \nu_i^{\kappa}
\beta(g)  \rangle
-
\sum_{\kappa\in\mathcal{T}^h}\int_{\kappa}
\langle v_i
\beta\, \partial_{x_i}\psi \rangle
 \\
 &=
\sum_{e\in\mathcal{S}^h}
\int_{e} \langle
\bm{v}\cdot\jump{\psi\,\hat{\beta}(g;v_{\nu})}\rangle
-
\sum_{\kappa\in\mathcal{T}^h}\int_{\kappa}\langle v_i
\beta\, \partial_{x_i} \psi \rangle
\end{split}
\end{equation}
where in the second equality $\beta(g)$ is replaced by any $\hat{\beta}(g;v_{\nu})$ in compliance with the consistency condition:
\begin{equation}
\label{eq:consistency}
\jump{\beta(g)}=0\quad\Rightarrow\quad\hat{\beta}(g;v_{\nu})=\beta(g)
\end{equation}  
Implicit in the identity in~\EQ{DGform2} is the assumption that~$\beta$ is sufficiently smooth within the elements to permit integration by parts and define traces on $\partial\kappa$. 
The edge distribution $\hat{\beta}(g;v_{\nu})$ is defined edge-wise and on each edge~$e$
it depends on $g$ only via $g_{\pm}$, viz. the restrictions of the traces of $g|_{\kappa_{\pm}}$ to~$e$. 
The function $\bm{v}\cdot\jump{\psi\,\hat{\beta}(g;v_{\nu})}$ in the ultimate expression in~\EQ{DGform2} 
can be conceived of as an upwind-flux weighted by the jump in $\psi$. It is to be noted that the domain of both
the upwind-flux and the jump $\jump{\psi}$ is $e\times\IR^D$.
On boundary edges, the external distribution corresponds to exogenous data in accordance with boundary conditions. 
Hence, any stationary solution to (\ref{eq:GalLevMomCls}) that is sufficiently regular in the aforementioned sense
satisfies
\begin{equation}
\label{eq:DGform1}
a(g;\psi)
=
s(g;\psi) \quad\forall\psi\in V^{h,p}(\Omega,\IM)
\end{equation}
with
\begin{align}
a(g;\psi)
&=
\sum_{e\in\mathcal{S}^h}
\int_{e}\langle
\bm{v}\cdot\jump{\psi \hat{\beta}(g;v_{\nu})}\rangle -
\sum_{\kappa\in\mathcal{T}^h}\int_{\kappa}
\langle
v_i \beta(g)\, \partial_{x_i}\psi \rangle
\label{eq:aDG}
\\
\label{eq:CollDG}
s(g;\psi)
&=
\sum_{\kappa\in\mathcal{T}^h}\int_{\kappa}
\langle \psi \CC(\beta(g))\rangle
\end{align}
The DGFE discretization of (\ref{eq:GalLevMomCls}) is obtained by replacing $g$ in~\EQ{DGform1} by an approximation $g^{h,p}_{\IM}$ in $V^{h,p}(\Omega,\IM)$ according to:
\begin{multline}
\label{eq:DGform}
\text{\em Find } g^{h,p}_{\IM}\in{}V^{h,p}(\Omega,\IM):
\quad
a(g{}^{h,p}_{\IM};\psi)
=
s(g{}^{h,p}_{\IM};\psi)
\quad\forall\psi\in{}V^{h,p}(\Omega,\IM)
\end{multline}

The edge distributions $\hat{\beta}$ in~\EQ{aDG} must be constructed such that the consistency condition~\EQ{consistency} 
holds and that the formulation~\EQ{DGform} is stable in some appropriate sense. We select the upwind edge 
distribution~\cite{Abdelmalik2016b}:
\begin{equation}
\label{eq:F+-}
\hat{\beta}(g;v_{\nu})
=
\begin{cases}
\beta(g{}_+)\quad&\text{if }v_{\nu_+}>0
\\
\beta(g{}_-)\quad&\text{if }v_{\nu_-}>0
\end{cases}
\end{equation}
In \cite{Abdelmalik2016b} it was shown that for suitable collision operators, (\ref{eq:F+-}) leads to an entropy stable formulation.

\section{Goal-oriented a-posteriori error estimation}
\label{sec:ErrEst}
If interest is restricted to a particular functional of the solution of~\EQ{Boltzmann}, the combined hierarchical and Galerkin structure of (\ref{eq:DGform}) may be used to derive an estimate of the error in the approximation of the goal functional.
In this section we will derive a computable a-posteriori goal-oriented error estimate
for~(\ref{eq:DGform}). We first present a formulation of the linearization of the DGFE moment
system~(\ref{eq:DGform}), which then serves as a basis for a computable a-posteriori error estimate in dual-weighted-residual (DWR) form~\cite{becker1996}.

We restrict ourselves to estimation of the {\em modeling error\/} that is incurred by limiting the dimension of the moment approximation in velocity dependence. We take the vantage point that the finite-element approximation space
in position dependence $V^{h,p}(\Omega)$ in~\EQ{DGform} is fixed and that~$\IM$ belongs 
to a nested sequence of moment spaces $\IE\subseteq\IM_0\subset\IM_1\subset\cdots\subset\IV$, where $\IV$
corresponds to a suitable closed normed vector space, and such that the 
sequence $\IM_k$ is asymptotically dense in~$\IV$, i.e. for all $g\in\IV$ and all $\epsilon>0$ there exists
a number $k_{\epsilon}\in\mathbb{Z}_{\geq0}$ such that $\inf_{g_k\in\IM_{k}}\|g-g_k\|_{\IV}<\epsilon$ for all $k\geq{}k_{\epsilon}$. Assuming that~\EQ{DGform} is well posed if~$\IM$ is replaced by~$\IV$ and denoting the corresponding 
solution by $g^{h,p}_{\IV}$, we are concerned with the error $J(g^{h,p}_{\IM})-J(g^{h,p}_{\IV})$ 
in the value of a goal functional $J:V^{h,p}(\Omega,{\IV})\to\IR$ evaluated at the approximation~$g^{h,p}_{\IM}$ 
according to~\EQ{DGform} for some finite-dimensional moment space~$\IM\subset\IV$.

The considered error estimate is based on linearization of~\EQ{DGform} at the approximation $g^{h,p}_{\IM}$.
We denote by
\begin{align*}
a'&:V^{h,p}(\Omega,\mathscr{V})\to(V^{h,p}(\Omega,\mathscr{V})\times{}V^{h,p}(\Omega,\mathscr{V}))^*
\\
s'&:V^{h,p}(\Omega,\mathscr{V})\to(V^{h,p}(\Omega,\mathscr{V})\times{}V^{h,p}(\Omega,\mathscr{V}))^*
\end{align*}
the Fr\'{e}chet derivatives of the semi-linear forms $a$ and $s$, respectively, with $(\cdot)^*$ the topological dual
of~$(\cdot)$. In particular, it holds that:
\begin{align}
a'[g](\delta g;\psi)
&=
\sum_{e\in\mathcal{S}^h}
\int_{e}\langle
\bm{v}\cdot\jump{\psi\hat{\beta}'[g](\delta g;v_{\nu})}\rangle
 -
\sum_{\kappa\in\mathcal{T}^h}\int_{\kappa}
\langle v_i
\beta'(g) \delta g\, \partial_{x_i}\psi \rangle
\label{eq:linaDG}
\\
\label{eq:linCollDG}
s'[g](\delta g;\psi)
&=
\sum_{\kappa\in\mathcal{T}^h}\int_{\kappa}
\langle \psi \CC'(\beta(g))\beta'(g)\delta g\rangle
\end{align}
where $\hat{\beta}'[g](\,\cdot\,;v_{\nu})$ represents the Fr\'{e}chet derivative of $\hat{\beta}(g;v_{\nu})$:
\begin{align*}
  \hat{\beta}'[g](\delta g;v_{\nu})
  =
  \begin{cases}
  \beta'(g{}_+) \delta g_+\quad&\text{if }v_{\nu^+}>0
  \\
  \beta'(g{}_-) \delta g_-\quad&\text{if }v_{\nu^-}>0
  \end{cases}
\end{align*}
Considering an approximation $g^{h,p}_{\IM}\in{}V^{h,p}(\Omega,\IM)$ to~$g^{h,p}_{\IV}$,
the error $\delta{}g=g^{h,p}_{\IV}-g^{h,p}_{\IM}$ satisfies:
\begin{multline}
\label{eq:linDGform}
a'[g^{h,p}_{\IM}](\delta g,\psi)
-
s'[g^{h,p}_{\IM}](\delta g,\psi)
=
-\operatorname{Res}[g^{h,p}_{\IM}](\psi)+o\big(\|\delta{}g\|_{L^{\infty}(\Omega,\IV)}\big)
\\
\forall{}\psi\in{}V^{h,p}(\Omega,\IV)
\end{multline}
as $\|\delta{}g\|_{L^{\infty}(\Omega,\IV)}\to{}0$, where $\operatorname{Res}:V^{h,p}(\Omega,\IV)\to(V^{h,p}(\Omega,\IV))^*$
is the residual functional according to:
\begin{equation}
\label{eq:ResDef}
\operatorname{Res}[g](\psi)=a(g;\psi)-s(g;\psi)
\end{equation}
and $\|\delta{}g\|_{L^{\infty}(\Omega,\IV)}=\sup\{\|\delta{}g(\bm{x},\cdot)\|_{\IV}:\bm{x}\in\Omega\}$.

We denote the Fr\'echet derivative of the target functional under consideration by 
$J':V^{h,p}(\Omega,\IV)\to(V^{h,p}(\Omega,\IV))^*$. In particular, if $J(\cdot)$ is given
by
\begin{equation}
\label{eq:Target}
  J(g) = \int_{\Omega} \langle \jmath_{\Omega}(g)\rangle + 
  \int_{\partial\Omega}\int_{v_n>0} \jmath_{\partial\Omega}(g)
\end{equation}
where $\jmath_{\Omega}:\IR\to\IR$ and $\jmath_{\partial\Omega}:\IR\to\IR$ are (possibly nonlinear) functions then
\begin{equation}
\label{eq:LinTarget}
J'[g](\delta g) = \int_{\Omega} \langle \jmath_{\Omega}'(g) \, \delta g \rangle 
+ 
\int_{\partial\Omega} \int_{v_n>0} \jmath_{\partial\Omega}'(g)\, \delta g
\end{equation}
To derive an estimate of the error in the target functional associated with the approximation $g^{h,p}_{\IM}$, 
we introduce the linearized {\em dual} or {\em adjoint} problem:
\begin{multline}
\label{eq:DGFEMAdj}
\text{\it Find }  
z\in{}V^{h,p}(\Omega,\IV): 
\\
a'[g^{h,p}_{\IM}](\delta g;z) - s'[g^{h,p}_{\IM}](\delta g;z) 
= 
J'[g^{h,p}_{\IM}](\delta g) \qquad \forall \delta g\in{}V^{h,p}(\Omega,\IV)
\end{multline}
The dual solution $z$ in~(\ref{eq:DGFEMAdj}) serves to construct an estimate of the error in the goal functional
in dual-weighted-residual form~\cite{becker1996} according to:
\begin{equation}
\label{eq:ApproxErrorRep}
\begin{aligned}
J(g^{h,p}_{\IM}) - J(g^{h,p}_{\IV}) 
&=-J'[g^{h,p}_{\IM}](\delta{}g)+o\big(\|\delta{}g\|_{L^{\infty}(\Omega,\IV)}\big)
\\ 
&=s'[g^{h,p}_{\IM}](\delta{}g;z) - a'[g^{h,p}_{\IM}](\delta{}g;z)+o\big(\|\delta{}g\|_{L^{\infty}(\Omega,\IV)}\big) 
\\
&=\operatorname{Res}[g^{h,p}_{\IM}](z)+o\big(\|\delta{}g\|_{L^{\infty}(\Omega,\IV)}\big) 
\end{aligned}
\end{equation}
as $\|\delta{}g\|_{L^{\infty}(\Omega,\IV)}\to{}0$.
The second identity in~\EQ{ApproxErrorRep} follows from~\EQ{DGFEMAdj}. The third identity follows from~\EQ{linDGform}.
The DWR error estimate is obtained by ignoring the $o(\|\delta{}g\|_{L^{\infty}(\Omega,\IV)})$ terms in the final expression in~\EQ{ApproxErrorRep}.

To elucidate the error estimate according to~\EQ{ApproxErrorRep}, we note that~\EQ{DGFEMAdj} can be regarded as an
approximation to the mean-value linearized dual problem:
\begin{multline}
\label{eq:NonLinAdj}
\text{\it Find } \bar{z}\in{}V^{h,p}(\Omega,\IV):
\\
\bar{a}(g^{h,p}_{\IM},g^{h,p}_{\IV};\delta g,\bar{z})
-
\bar{s}(g^{h,p}_{\IM},g^{h,p}_{\IV};\delta g,\bar{z})
= 
\bar{J}(g^{h,p}_{\IM},g^{h,p}_{\IV};\delta g)
\\
\forall \delta g\in{}V^{h,p}(\Omega,\IV):
\end{multline}
with
\begin{equation}
\begin{aligned}
\bar{a}(g^{h,p}_{\IM},g^{h,p}_{\IV};\delta g,\bar{z})
&=\int_0^1 a'[\theta g^{h,p}_{\IV}+(1-\theta)g^{h,p}_{\IM}](\delta g,\bar{z})\,d\theta
\\
\bar{s}(g^{h,p}_{\IM},g^{h,p}_{\IV};\delta g,\bar{z})
&=\int_0^1 s'[\theta g^{h,p}_{\IV}+(1-\theta)g^{h,p}_{\IM}](\delta g,\bar{z})\,d\theta
\\
\bar{J}(g^{h,p}_{\IM},g^{h,p}_{\IV};\delta g)
&=\int_0^1 J'[\theta g^{h,p}_{\IV}+(1-\theta)g^{h,p}_{\IM}](\delta g)\,d\theta
\end{aligned}
\end{equation}
For the mean-value linearized dual solution according to~\EQ{NonLinAdj}, the following (exact) error representation holds:
\begin{equation}
\label{eq:ErrRep}
J(g^{h,p}_{\IM}) - J(g^{h,p}_{\IV}) 
=
\operatorname{Res}[g^{h,p}_{\IM}](\bar{z})
\end{equation}
However, the mean-value linearized dual problem~\EQ{NonLinAdj} depends on $g^{h,p}_{\IV}$ and, accordingly, 
Equation~\EQ{ErrRep} does not provide a computable a-posteriori estimate. In the error estimate~\EQ{ApproxErrorRep}, the
mean-value linearized dual problem has been replaced by the linearized dual problem~\EQ{DGFEMAdj}, at the expense
of a linearization error $o(\|\delta{}g\|_{L^{\infty}(\Omega,\IV)})$ as~$\|\delta{}g\|_{L^{\infty}(\Omega,\IV)}\to0$. 

In practice, the dual problem~\EQ{DGFEMAdj} cannot be solved exactly and must again be approximated by a finite-element/moment
approximation. By the Galerkin orthogonality property of $g^{h,p}_{\IM}$ in~\EQ{DGform}, 
it holds that~$\operatorname{Res}[g^{h,p}_{\IM}](\psi)$ vanishes for all~$\psi\in{}V^{h,p}(\Omega,\IM)$. Hence,
for the dual solution in~\EQ{DGFEMAdj} an approximation space $V^{h,p}(\Omega,\IM_*)$ must be selected such that
$\IM_*\supset\IM$. More precisely, in the actual error estimate, the dual solution~$z$ is replaced by an approximation
$z^{h,p}_{\IM_*}$ according to:
\begin{multline}
\label{eq:DGFEMAdjh}
\text{\it Find }  
z^{h,p}_{\IM_*}\in{}V^{h,p}(\Omega,\IM_*): 
\\
a'[g^{h,p}_{\IM}](\delta g;z^{h,p}_{\IM_*}) - s'[g^{h,p}_{\IM}](\delta g;z^{h,p}_{\IM_*}) 
= 
J'[g^{h,p}_{\IM}](\delta g) 
\\
\forall \delta g\in{}V^{h,p}(\Omega,\IM_*)
\end{multline}
Typically, if the moment space $\IM$ is composed of polynomials up to order~$M$, then the refined 
space $\IM_*$ is selected such that it comprises all polynomials up to order~$M+1$ or~$M+2$.

\begin{remark}
The linearization that underlies the linearized dual problem~\EQ{DGFEMAdjh} can also serve in a 
Newton procedure to solve the nonlinear primal problem~\EQ{DGform}. 
\end{remark}

\begin{remark}
It is important to mention that the linearization of the semi-linear form~\EQ{aDG}
is significantly facilitated by basing the 
numerical flux on the upwind distribution according to~\EQ{F+-}. 
Alternatively, a discontinuous Galerkin approximation of~\EQ{GalLevMomCls} can be constructed
by first evaluating the velocity integrals and then introducing a DGFE approximation of the
resulting symmetric hyperbolic moment system. Such a DGFE formulation must then be equipped
with a numerical flux function (or {\em approximate Riemann solver\/}), e.g. according to Godunov's 
scheme~\cite{Godunov1959}, Roe's scheme~\cite{Roe:1981fj} or Osher's scheme~\cite{Osher:1982tg}. 
However, these numerical flux functions generally depend in an intricate manner on the left and right
states via the eigenvalues and eigenvectors of the flux Jacobian, Riemann invariants, etc., 
which impedes differentiation of the resulting semi-linear form. Determing the derivative of
the upwind distribution in~\EQ{F+-} and, in turn, of the semi-linear form~\EQ{aDG} is a straightforward 
operation.
\end{remark}

\begin{remark}
Implicit to the error representation in~(\ref{eq:ApproxErrorRep}) is the assumption that 
the nonlinear primal problem~(\ref{eq:DGform}) and the linearized dual problem~(\ref{eq:DGFEMAdj})
are well posed. A rigorous justification of this assumption is technical and beyond the scope of this work.
Let us mention however that~\EQ{GalLevMomCls} represents a symmetric hyperbolic systems and, hence, provided
with suitable auxiliary conditions it is linearly well posed for sufficiently smooth solutions~\cite{Majda:1984wj}. In particular, the linearized system corresponds to a Friedrichs system~\cite{Friedrichs1958}. Accordingly, 
the linearization of~\EQ{DGform} corresponds to a DGFE approximation of a Friedrichs system equipped with a 
proper upwind flux. Discontinuous Galerkin approximations of Friedrichs systems are generally well posed;
see, for instance, \cite{Ern:2006kl,Ern:2006pi,Ern:2008vl}. Well-posedness of  the linearized adjoint problem (\ref{eq:DGFEMAdjh}) follows directly from well-posedness of the corresponding linearized primal problem;
see~\cite[proposition A.2]{melenk1999}.
\end{remark}

\section{Goal-oriented adaptive algorithm}
\label{sec:AdaptAlg}
The computable error estimate (\ref{eq:ApproxErrorRep}) can be used to direct an adaptive algorithm following the
standard SEMR (\texttt{Solve} $\to$ \texttt{Estimate} $\to$ \texttt{Mark} $\to$ \texttt{Refine}) process; see
for instance~\cite{dorfler1996,Brummelen:2017rr,Nochetto:2012hl}. The marking step comprises
a decomposition of the error estimate~(\ref{eq:ApproxErrorRep}) into local element-wise contributions,
and a subsequent marking of elements that provide the dominant contributions to the error. To enhance
the efficiency of the adaptive algorithm, we consider a marking strategy that accounts for cancellation effects.
The refinement process consists in locally raising the number of moments in the elements that have been
marked for refinement. By repeated application of SEMR, the adaptive algorithm aims to adapt the number of moments 
locally in each element to obtain an optimal approximation to the quantity of interest. 

By virtue of the local nature of the discontinous Galerkin approximation in position dependence, the element-wise
decomposition of the error estimate according to the ultimate expression in~(\ref{eq:ApproxErrorRep}) is straightforward.
Denoting by $\{\Lambda_{\kappa,i}(\bm{x},\bm{v})\}$ a basis of the approximation space $V^{h,p}(\Omega,\IM_*)$ for the dual solution such that the support of each function $\Lambda_{\kappa,i}$ is confined to the 
element~$\kappa\in\mathcal{T}^h$, it holds that
\begin{equation}
\label{eq:zetadef}
\operatorname{Res}[g^{h,p}_{\IM}](z^{h,p}_{\IM_*})
=
\sum_{\kappa\in\mathcal{T}^h}
\underbrace{\sum_{i\in\mathcal{I}_{\kappa}}
\operatorname{Res}[g^{h,p}_{\IM}](\Lambda_{\kappa,i})\,\sigma_{\kappa,i}}_{\zeta_{\kappa}}
\end{equation}
where $\mathcal{I}_{\kappa}$ is an index set corresponding to element~$\kappa\in\mathcal{T}^h$ 
and $\sigma_{\kappa,i}$ are the weights of $z^{h,p}_{\IM_*}$ relative to~$\{\Lambda_{\kappa,i}(\bm{x},\bm{v})\}$.
Indeed, the error contributions $\{\zeta_{\kappa}\}$ are directly associated with the elements.

Conventionally, to mark elements for refinement an upper bound for the error estimate (\ref{eq:ApproxErrorRep}) is constructed based on the absolute value of $\zeta_{\kappa}$ and the triangle inequality:
\begin{equation}
\label{eq:Triangle}
\Big|\operatorname{Res}[g^{h,p}_{\IM}](z^{h,p}_{\IM_*})\Big|
\leq
\sum_{\kappa\in\mathcal{T}^h}|\zeta_{\kappa}|
\end{equation}
see, for instance,~\cite{hartmann2002,hartmann2003,houston2001}. Elements are then marked for refinement
according to their contribution  $|\zeta_{\kappa}|$ to the upper bound, e.g., 
following the D\"{o}rfler marking strategy~\cite{dorfler1996}, which selects a minimal set of 
elements $\mathcal{T}^h_{*}\subset\mathcal{T}^h$ such that:
\begin{equation}
\label{eq:Dorfler}
\sum_{\kappa\in \mathcal{T}^h_{*}} |\zeta_{\kappa}| \geq c\sum_{\kappa\in\mathcal{T}^h} |\zeta_{\kappa}|
\end{equation}
for some $c\in(0,1]$. However, previous work for first order hyperbolic systems~\cite{hartmann2002,hartmann2003} 
suggests that the upper bound provided by the triangle inequality may not be sharp due to the loss of inter-element cancellations. In this work we aim to exploit such cancellation errors. To that end, we decompose $\mathcal{T}^h$ according to the sign of the local error indicators $\zeta_{\kappa}$ relative to the error estimate $\sum_{\kappa\in\mathcal{T}^h} \zeta_{\kappa}$ according to~$\mathcal{T}^h = \mathcal{T}^h_+ \cup \mathcal{T}_-^h$ where
\begin{align}
\begin{split}
	\mathcal{T}_+^h &= \bigg\{\kappa\in\mathcal{T}^h: \text{sgn}\left(\zeta_{\kappa}\right) = \text{sgn}\bigg(\sum_{\kappa\in\mathcal{T}^h} \zeta_{\kappa}\bigg) \bigg\}\\
	\mathcal{T}_-^h &= \bigg\{\kappa\in\mathcal{T}^h: \text{sgn}\left(\zeta_{\kappa}\right) \neq \text{sgn}\bigg(\sum_{\kappa\in\mathcal{T}^h} \zeta_{\kappa}\bigg) \bigg\}
\end{split}
\end{align} 
That is, the elements in~$\mathcal{T}_+^h$ (resp. $\mathcal{T}_-^h$) are those whose local error contribution increases (resp. decreases) $|\sum_{\kappa\in\mathcal{T}^h}\zeta_{\kappa}|$. We propose to mark a minimal set of elements 
$\mathcal{T}^h_*\subset\mathcal{T}_+^h$ such that
\begin{equation}
\label{eq:mark}
\bigg|\sum_{\kappa\in \mathcal{T}_*^h} \zeta_{\kappa}\bigg| \geq c\,\bigg|\sum_{\kappa\in\mathcal{T}^h} \zeta_{\kappa}\bigg|
\end{equation}
for some $c\in(0,1]$. To elucidate the manner in which this approach accounts for inter-element cancellations, 
we note that the proposed marking strategy in~(\ref{eq:mark}) may be equivalently understood as selecting elements 
in~$\mathcal{T}^h_+$ with the largest (in magnitude) error contributions from the complement of the 
the largest (in cardinality) subset of elements in~$\mathcal{T}^h_+$ of which the contribution can be 
cancelled against the aggregated contribution of elements in~$\mathcal{T}_{-}^h$. More precisely, denoting by
$\widetilde{\mathcal{T}}{}^h_+\subset\mathcal{T}^h_+$ a maximal set of elements such that
\begin{equation}
\label{eq:mark'}
\Big|\sum_{\kappa\in\widetilde{\mathcal{T}}{}_{+}^h} \zeta_{\kappa}\Big|\leq\Big|\sum_{\kappa\in\mathcal{T}_-^h} \zeta_{\kappa}\Big|
\end{equation}
it holds that
\begin{multline}
\label{eq:Bnd}
\Big|\operatorname{Res}[g^{h,p}_{\IM}](z^{h,p}_{\IM_*})\Big|
=
\bigg|
\sum_{\kappa\in\widetilde{\mathcal{T}}{}_{+}^h}\zeta_{\kappa} 
+
\sum_{\kappa\in\mathcal{T}^h_+\setminus\widetilde{\mathcal{T}}{}_{+}^h}\zeta_{\kappa} 
+
\sum_{\kappa\in\mathcal{T}^h_-}\zeta_{\kappa} 
\bigg|
\\
=
\sum_{\kappa\in\widetilde{\mathcal{T}}{}_{+}^h}|\zeta_{\kappa}| 
+
\sum_{\kappa\in\mathcal{T}^h_+\setminus\widetilde{\mathcal{T}}{}_{+}^h}|\zeta_{\kappa}| 
-
\sum_{\kappa\in\mathcal{T}^h_-}|\zeta_{\kappa}| 
\leq
\sum_{\kappa\in\mathcal{T}^h_+\setminus\widetilde{\mathcal{T}}{}_{+}^h}|\zeta_{\kappa}| 
\end{multline}
Hence, the aggregated error contributions of the elements in~$\mathcal{T}^h_+\setminus\widetilde{\mathcal{T}}{}_{+}^h$
yield an upper bound to the error estimate. The set of marked elements~$\mathcal{T}^h_*$ corresponds
to the smallest subset (in cardinality) of~$\mathcal{T}^h_+\setminus\widetilde{\mathcal{T}}{}_{+}^h$ such that~\EQ{mark} holds.
Hence, one can conceive of the proposed marking strategy in~\EQ{mark} as a procedure that first accounts for 
inter-element cancellations and subsequently applies a D\"orfler-type marking to the remaining elements.

The SEMR algorithm based on~\EQ{mark} is summarized in Algorithm~\ref{alg:SEMR}. One first defines a sequence
of moment spaces $\{\IM_r\}_{r\in\mathbb{Z}_{\geq0}}$ and $\{\IM_{r*}\}_{r\in\mathbb{Z}_{\geq0}}$ for the primal and dual problems, respectively. It must
hold that $\IM_{r*}\supset\IM_r$ to avoid that the error estimate vanishes due to Galerkin orthogonality; see 
Section~\SEC{ErrEst}. Next, the element-wise hierarchical rank of the moment-system approximation is initialized
at the basic level $r_{\kappa}=0$. In the iterative process, one first constructs the (possibly non-uniform)
approximation spaces for the primal and dual problem,~viz.
\begin{equation}
V^{h,p}\big(\Omega,\{\IM_{r_{\kappa}}\}\big)
=
\big\{g\in{}V^{h,p}(\Omega,\IV):g|_{\kappa}\in{}\mathcal{P}_p(\kappa,\IM_{r_{\kappa}})\big\}
\end{equation}
and likewise for the dual problem, and then solves the nonlinear primal problem~\EQ{DGform} and 
the linearized dual problem~\EQ{DGFEMAdjh} for the approximate primal and dual solutions. Based on
the approximate dual solution and the residual corresponding to the approximate primal solution,
the error contributions $\zeta_{\kappa}$ can be computed and the error estimate $\texttt{est}$ 
can be assembled. If the estimate satisfies the prescribed tolerance, then the algorithm terminates.
Otherwise, the algorithm proceeds by marking a minimal set of elements $\mathcal{T}^h_*$ following
the D\"orfler marking with cancellations in~\EQ{mark}. In these marked elements, the approximation
is refined by incrementing the hierarchical rank of the moment approximation.
\begin{Algorithm}
\begin{center}
\begin{boxedminipage}{\textwidth}
\noindent
{\small
\begin{algorithmic}[1]
\State define $\IM_0,\IM_1,\IM_2,\ldots$
\Comment{sequence of primal moment spaces}
\State define $\IM_{0*}\supset\IM_0,\IM_{1*}\supset\IM_,\IM_{2*}\supset\IM_2,\ldots$
\Comment{sequence of dual moment spaces}
\State set $r_{\kappa}=0$ for all $\kappa\in\mathcal{T}^h$ 
\Comment{initialize element-wise moment rank}
\State set $\texttt{tol}>0$
\Comment{define tolerance}
\Loop
\State solve nonlinear primal problem~\EQ{DGform} for $g^{h,p}_{\IM}\in{}V^{h,p}\big(\Omega,\{\IM_{r_{\kappa}}\}\big)$
\Comment{solve}
\State solve linearized dual problem~\EQ{DGFEMAdjh} for $z^{h,p}_{\IM_*}\in{}V^{h,p}\big(\Omega,\{\IM_{r_{\kappa}*}\}\big)$
\State determine element-wise error indicators $\zeta_{\kappa}$ according to~\EQ{zetadef}
\State determine error estimate ${\texttt{est}}=\sum_{\kappa\in\mathcal{T}^h}\zeta_{\kappa}$
\Comment{estimate}
\If{$\texttt{est}<\texttt{tol}$}
\State break
\Else
\State mark a minimal subset of elements $\mathcal{T}{}^h_*$ according to~\EQ{mark}
\Comment{mark}
\For{$\kappa\in\mathcal{T}{}^h_*$}
\State $r_{\kappa} \gets r_{\kappa}+1$
\Comment{refine}
\EndFor
\EndIf
\EndLoop
\end{algorithmic}}
\end{boxedminipage}
\end{center}
\caption{The SEMR algorithm for goal-oriented element-wise model refinement of the moment-system approximation.\label{alg:SEMR}}
\end{Algorithm}

It is noteworthy that the adaptive algorithm admits a reinterpretation as a type\nobreakdash-A heterogeneous multiscale 
method (HMM)~\cite{e2003}. Multiscale methods of type~A introduce a decomposition of the spatial domain into a region where a microscopic description based on a sophisticated model is required, and a complementary
region where a macroscopic description based on a simple model suffices. The adaptive strategy in Algorithm~\ref{alg:SEMR}
forms an automatic partitioning of the domain~$\Omega$ into subregions (corresponding to collections of elements) on which
models of different levels of sophistication are applied, viz. the moment systems associated with the different hierarchical
ranks of the moment spaces~$\IM_0,\IM_1,\IM_2,\ldots$ Instead of two models, i.e. one microscopic model and one macroscopic model, the moment systems furnish a sequential hierarchy of models. The SEMR procedure in Algorithm~\ref{alg:SEMR} fully automates the domain decomposition strategy and model selection in such a manner that a (quasi\nobreakdash-)optimal approximation of the quantity of interest is obtained.

If the baseline moment space~$\IM_0$ coincides with 
the collision invariants~$\IE$ and the background distribution~$\mathcal{B}$ is appropriately chosen, 
then the lowest rank model in the sequence corresponds to the Euler equations.
On the other hand, the approximation provided by the sequence of moment-system approximations converges 
(formally) to the solution of the underlying Boltzmann equation as the hierarchical rank is refined, i.e. as $r\to\infty$. 
The goal-adaptive DGFE moment approximation hence implements a type\nobreakdash-A hierarchical HMM for combining
the Euler equations (macroscale model) and the Boltzmann equation (microscale model), as suggested 
in~\cite[Sec.~2.3.2]{e2003}.

\section{Numerical Results}
\label{sec:NumRes}
To illustrate the properties of the proposed goal-oriented model-adaptive strategy in section \ref{sec:AdaptAlg} for the discontinuous Galerkin finite-element moment method~\EQ{DGform}, we present numerical experiments for heat transfer and 
shock-structure problems in one dimension; see, for example,~\cite{gamba2010}. 

The problem specification and moment-system approximation must be completed first by specifying the collision operator and closure relation. We restrict ourselves here to the standard BGK collision operator~\cite{Bhatnagar:1954hc}, viz.
\begin{align}\label{eq:BGK}
  \CC(f) = \tau^{-1}(\mathcal{M}_f-f)
\end{align}
where $\mathcal{M}_f$ denotes the local equilibrium Maxwellian (\ref{eq:Maxwellian}) having the same invariant moments as~$f$ and~$\tau^{-1}$ is a relaxation rate. We adopt the relaxation parameter in accordance with the hard-sphere collision process of Bird~\cite{Bird:1994} 
\begin{equation}
\tau = (5\lambda/16)(2\pi\rho/p)^{1/2}
\end{equation}
with $\lambda$ the mean free path.
We consider discontinuous Galerkin finite-element approximation spaces of polynomial degree $p = 0$, i.e. element-wise constant approximations in position dependence. To solve the DGFE approximation (\ref{eq:DGform}) we use a Newton procedure based on the linearized DGFE approximation in (\ref{eq:linDGform}). To illustrate the usefulness of the adaptive procedure in capturing non-equilibrium flow phenomena, in the sequel, the goal-oriented adaptive algorithm considers the average of the heat flux 
\begin{equation}
\label{eq:SSGoal}
J(g)=\int_{\Omega}\big\langle (v-U)^3 \beta(g(x,v)) \big\rangle\, dx,
\end{equation}
as the quantity of interest. In all cases we choose a greedy refinement strategy and set the refinement fraction 
in~\EQ{mark} to~$c=1$, i.e. all elements that remain after cancellations have been accounted for are refined.

\subsection{Heat Transfer Problem}
The first test case pertains to the so-called heat transfer problem \cite{gamba2010}. This test case is set on a
unit interval~$\Omega=(0,1)$. We consider full accommodation boundary conditions at the left and right boundaries of the domain, with boundary data corresponding to uniform Maxwellian distributions with different temperatures:
\begin{equation}
\label{eq:HTdata}
  f_{\textsc{ht}}(v)=
  \begin{cases}
    \mathcal{M}_{(\rho^{\textsc{ht}}_l,0,\theta^{\textsc{ht}}_l)}(v) &\qquad x=0,\ v>0 \\
    \mathcal{M}_{(\rho^{\textsc{ht}}_r,0,\theta^{\textsc{ht}}_r)}(v) &\qquad x=1,\ v<0
  \end{cases}
\end{equation}
The left and right boundary densities, $\rho^{\textsc{ht}}_l$ and $\rho^{\textsc{ht}}_r$, respectively, are 
determined from the mass impermeability condition:
\begin{equation}
\label{eq:imperm}
\int_{v_n<0}v_nf_{\textsc{ht}}(v)\,dv+\int_{v_n>0}v_n\beta(g(v))\,dv = 0
\end{equation}
Condition (\ref{eq:imperm}) imposes that the entering and exiting mass fluxes on~$\partial\Omega$ cancel.

For this problem we consider the renormalization map (\ref{eq:phimap}) with~$N=1$. Note that for any $N>0$, the resulting moment-system is non-linear due to the non-linearity of the collision operator (\ref{eq:BGK}) and the non-linearity of the renormalization map (\ref{eq:phimap}). We consider a spatially non-uniform background distribution $\mathcal{B}^{\textsc{ht}}(x,v)=\mathcal{M}_{(\rho^{\textsc{ht}}(x),0,\theta^{\textsc{ht}}(x))}(v)$ where 
\begin{align}
\begin{split}
\theta^{\textsc{ht}}(x) &= \theta^{\textsc{ht}}_l+(\theta_r^{\textsc{ht}}-\theta^{\textsc{ht}}_l)x,\\
\rho^{\textsc{ht}}(x)&=(\theta^{\textsc{ht}}_l+\theta^{\textsc{ht}}_r)/(2\theta^{\textsc{ht}}(x))
\end{split}
\end{align} 
It is noteworthy that $\rho^{\textsc{ht}}(x)$ and $\theta^{\textsc{ht}}(x)$ can be conceived of as continuum approximations 
satisfying~\cite{bobylev1996}:
\begin{equation}
	\theta^{\textsc{ht}}(0)=\theta^{\textsc{ht}}_l, \quad \theta^{\textsc{ht}}(1)=\theta_r^{\textsc{ht}} \quad\text{and}\quad \frac{d}{dx}\frac{1}{\rho^{\textsc{ht}}}\frac{d\theta^{\textsc{ht}}}{dx}=C
\end{equation}
for some constant~$C$.
The adaptive algorithm is initiated with a uniform moment approximation that includes moments up to 
order~$M_{\kappa}=4$. The dual
solution is approximated based on~\EQ{DGFEMAdjh} using a moment approximation that is refined by raising the order 
in each element to $M_{\kappa*}=M_{\kappa}+2$.
We consider the heat transfer problem with Knudsen number $\Kn=10^{-3}$ and $\theta^{\textsc{ht}}_l=1$ 
and~$\theta_r^{\textsc{ht}}=1.2\theta_l^{\textsc{ht}}$. The domain is covered with a uniform mesh with~$10^3$ elements. 

\begin{figure}
\centering
\includegraphics[width=\textwidth]{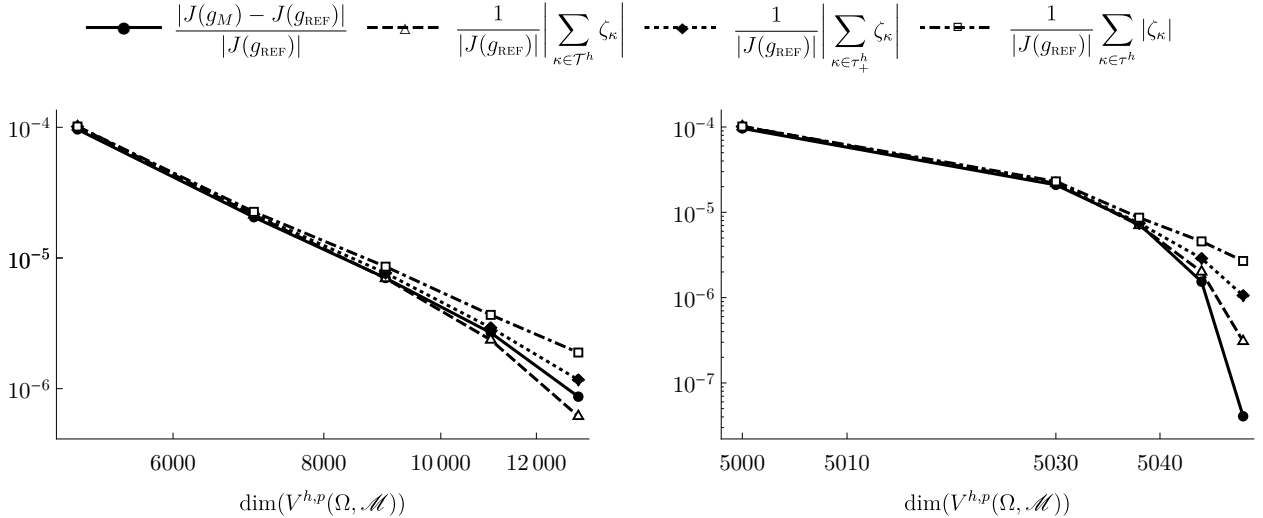}
\caption{
Error in the goal-functional according to~\EQ{SSGoal} ({\em solid\/}), goal-oriented error estimate ({\em dashed\/}),
error bound including cancellations ({\em dotted\/}) and conventional error bound ({\em dash-dot\/})
for the heat transfer problem versus the number of degree of freedom, for uniform  moment-order refinement (\textit{left}) 
and goal-adaptive refinement by Algorithm~\ref{alg:SEMR} (\textit{right}). 
\label{fig:HTadapt}}
\end{figure}

\begin{figure}
\centering
\includegraphics[width=\textwidth]{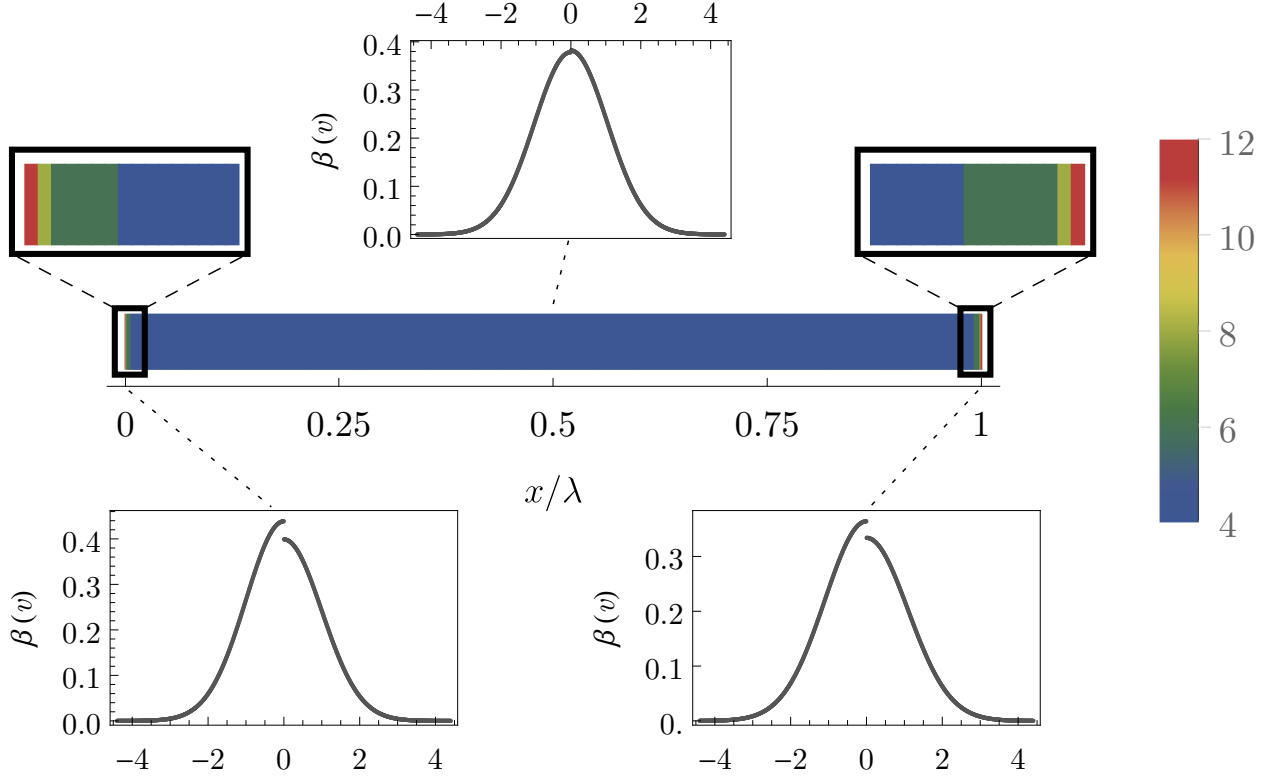}
\caption{Order of moment approximation in the computational domain for the goal-adaptive approximation of the solution to the heat transfer problem with zooms of the regions near the boundaries. The inserted panels 
display the upwind distribution according to~\EQ{F+-}.\label{fig:HTvdeg}}
\end{figure}

Figure~\ref{fig:HTadapt} presents the error $|J(g^{h,p}_{\IM})-J(g^{h,p}_{\textsc{ref}})|$ 
with respect to a reference result $J(g^{h,p}_{\textsc{ref}})$ based on a spatially uniform 
approximation with moments up to order $M_{\kappa}=14$. 
Figure~\ref{fig:HTadapt} also displays the error estimate, $|\sum_{\kappa}\zeta_{\kappa}|$, the upper bound
including cancellations according to the ultimate expression in~\EQ{Bnd}. Moreover, we plot 
the conventional upper bound of the error estimate based on the triangle inequality according to~(\ref{eq:Triangle}). The left panel plots the aforementioned error estimates and bounds versus the number of degrees of freedom 
under uniform refinement in the number of moments for $M_{\kappa}\in\{4,6,8,10,12\}$. The right panel displays the 
corresponding results for the goal-adaptive approximation obtained by means of Algorithm~\ref{alg:SEMR}. The results 
in Figure~\ref{fig:HTadapt} convey that under uniform refinement, $8\times{}10^3$ additional degrees of freedom 
are required to achieve a relative error of~$10^{-6}$ in the heat flux~\EQ{SSGoal}. The goal-adaptive refinement strategy 
only requires $50$ additional degrees of freedom to achieve the same accuracy. These results hence provide a clear indication
of the efficiency gain that can be obtained by the goal-oriented adaptive-refinement process. It is to noted that
the moment order in the adaptive approximation has locally reached $M_{\kappa}=12$ in the final approximation; see also Figure~\FIG{HTvdeg}.  

Comparison of the upper bounds in Figure~\ref{fig:HTadapt} confirms that the proposed upper bound 
in~(\ref{eq:Bnd}) is sharper than the standard triangle inequality~\EQ{Triangle}, which illustrates the effect
of cancellations. In particular, the deviation between the bounds becomes more pronounced as the approximation is
refined and the number of moments increases, both for the uniform approximation and for the goal-adaptive approximation.

To illustrate the spatial distribution of the moment refinement in the goal-adaptive approximation, 
Figure~\ref{fig:HTvdeg} plots the order of the moments in every element in the final step of the adaptive algorithm.
In addition, Figure \ref{fig:HTvdeg} displays the upwind distributions $\hat{\beta}(g^{h,p}_{\{\IM_{\kappa}\}};v_{\nu})$ 
on the left and right boundaries and at the element interface in the center of the spatial domain.
The results in Figure~\ref{fig:HTvdeg} show that most of the moment refinements occur in the regions near the boundaries. This can be attributed to the fact that the solution exhibits large jumps at the domain boundaries, as indicated by the plots of $\hat{\beta}(g^{h,p}_{\{\IM_{\kappa}\}};v_{\nu})$. These jumps represent non-equilibrium effects due
to incompatibility of the solution with the Maxwellian equilibrium distributions at the boundary.

\begin{figure}
\centering
\subfloat[Primal distribution $\beta(g_{\IM}^{h,p})$]{\includegraphics[width=0.8\textwidth]{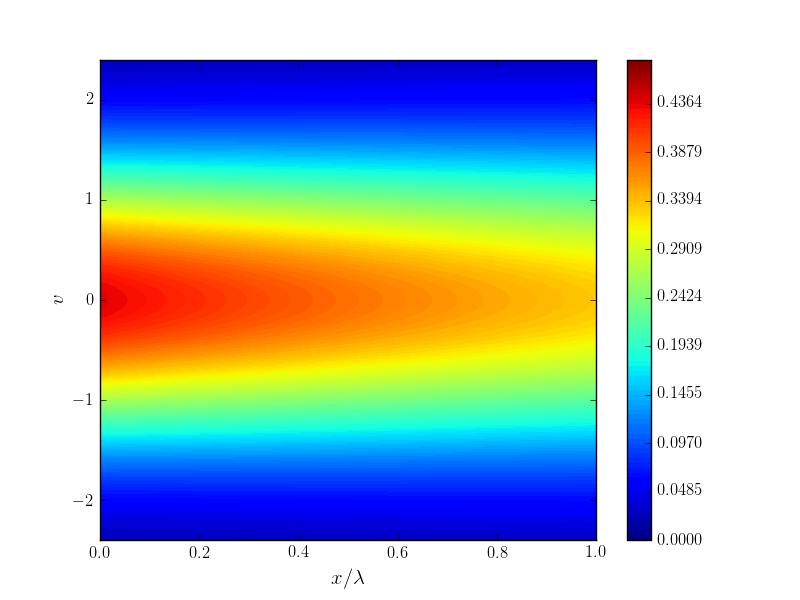}}\,
\subfloat[Dual solution $ z_{\IM_*}^{h,p}$]{\includegraphics[width=0.8\textwidth]{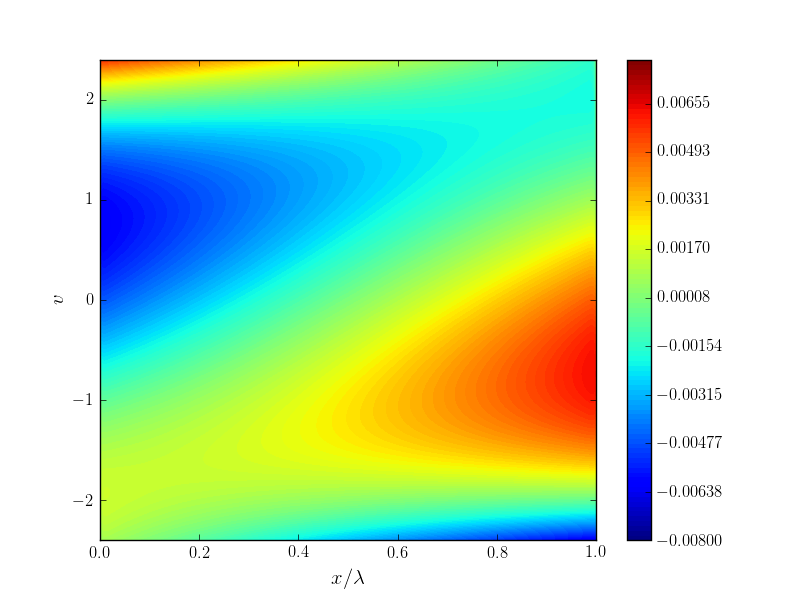}}
\caption{Goal-adaptive approximation of the primal distribution $\beta(g_{\{\IM_{\kappa}\}}^{h,p})$ 
for the heat transfer problem ({\it top}) and corresponding approximate dual solution $z^{h,p}_{\{\IM_{\kappa*}\}}$
({\it bottom}) in the final step of the adaptive algorithm\label{fig:HTPD}.}
\end{figure}

To further elucidate the refinement pattern in Figure~\ref{fig:HTvdeg}, Figure~\ref{fig:HTPD} displays
the approximation of the primal distribution $\beta(g^{h,p}_{\{\IM_{\kappa}\}})$ ({\em top\/}) and 
the dual distribution $z^{h,p}_{\{\IM_{\kappa*}\}}$ ({\em bottom\/}) in the final step of the adaptive algorithm.
The dual solution assigns most weight to the regions near the boundaries. In combination with the fact that large
residuals occur near the boundaries on account of incompatibility of the boundary data with the solution (see Figure~\ref{fig:HTvdeg}), the error contributions $\zeta_{\kappa}$ are mostly localized in the vicinity of the boundaries.

\subsection{Shock Structure Problem}
The second test case that we consider pertains to the so-called shock-structure problem~\cite{gamba2010} on a spatial domain $\Omega=(-40\lambda,40\lambda)$. This test case concerns a Riemann problem with boundary data corresponding to uniform Maxwellian distributions:
\begin{equation}
  f_{\textsc{ss}}=
  \begin{cases}
    \mathcal{M}_{(\rho^{\textsc{ss}}_l,u^{\textsc{ss}}_l,\theta^{\textsc{\textsc{ss}}}_l)}(v) &\qquad x=-40\lambda;\ v>0 \\
    \mathcal{M}_{(\rho^{\textsc{ss}}_r,u^{\textsc{ss}}_r,\theta^{\textsc{ss}}_r)}(v) &\qquad x=\phantom{-}40\lambda;\ v<0
  \end{cases}
\end{equation}
where the density, mean velocity and temperature on both sides of the shock are related by the Rankine-Hugoniot conditions \cite{gombosi1994}:
\begin{align}
\label{eq:ICdata}
\begin{split}
	\rho^{\textsc{ss}}_r &= \rho^{\textsc{ss}}_l\frac{(\gamma+1)\mathrm{Ma}^2}{2+(\gamma-1)\mathrm{Ma}^2}\\
	\theta^{\textsc{ss}}_r &= \theta^{\textsc{ss}}_l \frac{\rho_l}{\rho_r} \frac{2\gamma \mathrm{Ma}^2-(\gamma-1)}{\gamma+1}\\
	u_l^{\textsc{ss}} &= \mathrm{Ma} \sqrt{\gamma\theta_l^{\textsc{ss}}}\\
	u_r^{\textsc{ss}} &= u_l^{\textsc{ss}} \frac{2+(\gamma-1)\mathrm{Ma}^2}{(\gamma+1)\mathrm{Ma}^2}
\end{split}
\end{align}
with $\mathrm{Ma}$ denoting the Mach number and $\gamma=1+2/n$ the so-called adiabatic exponent for a perfect gas whose molecules have $n$ degrees freedom \cite{gombosi1994}. For this problem we consider the renormalization map (\ref{eq:phimap}) with $N=2$ and a background distribution $\mathcal{B}^{\textsc{ss}}(x,v)=\mathcal{M}_{(\rho^{\textsc{ss}}(x),u^{\textsc{ss}}(x),\theta^{\textsc{ss}}(x))}(v)$ that derives from the boundary data in (\ref{eq:ICdata}) as:
\begin{align}\label{eq:bimod}
\begin{split}
	\frac{\rho^{\textsc{ss}}(x)}{\rho_l^{\textsc{ss}}}&=X(x)+\frac{u_l^{\textsc{ss}}}{u_r^{\textsc{ss}}}\big(1-X(x)\big)=\frac{u_l^{\textsc{ss}}}{u^{\textsc{ss}}(x)}
	\\
	\frac{\theta^{\textsc{ss}}(x)}{\theta^{\textsc{ss}}_l}&=\frac{\rho^{\textsc{ss}}_l}{\rho^{\textsc{ss}}(x)}\left(X(x)+\frac{\theta_r^{\textsc{ss}}}{\theta_l^{\textsc{ss}}}\frac{u_l^{\textsc{ss}}}{u_r^{\textsc{ss}}}\big(1-X(x)\big)\right)
	\\ 
	&\phantom{=}+ \frac{\gamma}{3}\left(\frac{\rho^{\textsc{ss}}_l}{\rho^{\textsc{ss}}(x)}\right)^2\left(1-\frac{u_r^{\textsc{ss}}}{u_l^{\textsc{ss}}}\right)^2\mathrm{Ma}^2\big(1-X(x)\big)X(x)
\end{split}
\end{align}
with the interpolation function
\begin{equation}
\label{eq:X(x)}
X(x) = \frac{1}{2}-\frac{1}{2}\tanh\left(\frac{2x}{40\lambda}\right)
\end{equation}
Let us mention that the interpolation function in~\EQ{X(x)} has been chosen such that for $x=-40\lambda$ (resp. $x=40\lambda$)
it holds that $X(x)$ is close to~$0$ (resp. close to~$1$), but it is otherwise arbitrary.
The background distribution~$\mathcal{B}^{\textsc{ss}}(x,v)$ is understood as a local Maxwellian approximation of a distribution that interpolates the boundary data~(\ref{eq:ICdata}) using the interpolation function $X(x)$, similar to 
the so-called Mott-Smith approximation~\cite{anan1969,gombosi1994}.

The adaptive algorithm is initiated with a spatially uniform moment approximation of degree $M_{\kappa}=4$. 
The linearized dual problem is approximated using a moment approximation that is refined by locally raising the 
order to $M_{\kappa*}=M_{\kappa}+4$. In this case we opt to apply 
$M_{\kappa*}=M_{\kappa}+4$ instead of~$M_{\kappa*}=M_{\kappa}+2$ to improve the accuracy of the error estimate.
The dual solution exhibits non-smooth behavior near the boundaries and, accordingly, insufficient resolution 
in velocity dependence leads to an inferior error estimate. 

\begin{figure}
\centering
\includegraphics[width=\textwidth]{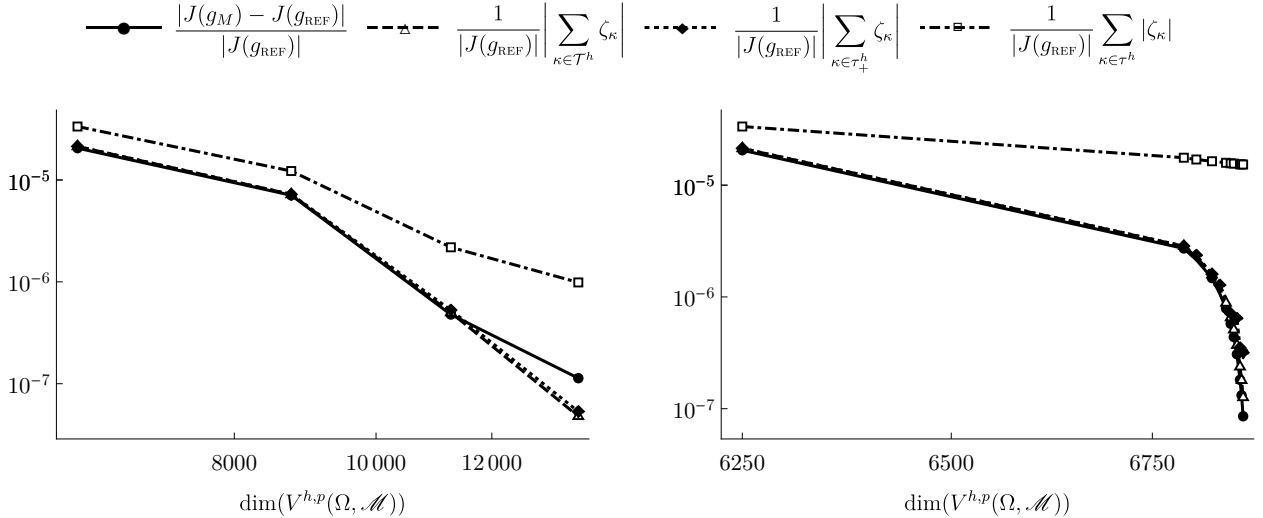}
\caption{
Error in the goal-functional according to~\EQ{SSGoal} ({\em solid\/}), goal-oriented error estimate ({\em dashed\/}),
error bound including cancellations ({\em dotted\/}) and conventional error bound ({\em dash-dot\/})
for the shock-structure problem versus the number of degree of freedom, for uniform  moment-order refinement (\textit{left}) and goal-adaptive refinement (\textit{right}). 
\label{fig:SSadapt}}
\end{figure}

We consider the shock structure problem with Mach number $\mathrm{Ma}=1.4$ and mean free path $\lambda=3.67\times{}10^{-3}$. The computational domain is covered with a uniform mesh of~$1250$ elements. 
Figure~\ref{fig:SSadapt} shows the error $|J(g^{h,p}_{\{\IM_{\kappa}\}})-J(g^{h,p}_{\textsc{ref}})|$ relative to the reference result $J(g^{h,p}_{\textsc{ref}})$ based on a spatially uniform approximation with moments up to order $M_{\kappa}=12$.
In addition, the figure displays the error estimate, the upper bound
including cancellations according to the ultimate expression in~\EQ{Bnd} and 
the conventional upper bound~(\ref{eq:Triangle}). The left panel presents the results for
uniform refinement in the number of moments for $M_{\kappa}\in\{4,6,8,10\}$. The right panel presents 
results for the goal-adaptive approximation.
Figure~\ref{fig:SSadapt} shows that for the shock structure problem, uniform refinement requires more than $7500$ 
additional degrees of freedom to reduce the relative error to $10^{-7}$. The adaptive strategy only requires~$616$ additional degrees of freedom to reach the same relative error.
The results reaffirm that significant gains in efficiency can be obtained by means of the goal-adaptive refinement strategy.
It may be noted that for the considered shock-structure test case, the conventional error bound 
derived from the triangle inequality~\EQ{Triangle} is very loose, while the bound~(\ref{eq:Bnd}) 
that accounts for cancellations is sharp relative to the error estimate.

\begin{figure}
\begin{center}
\includegraphics{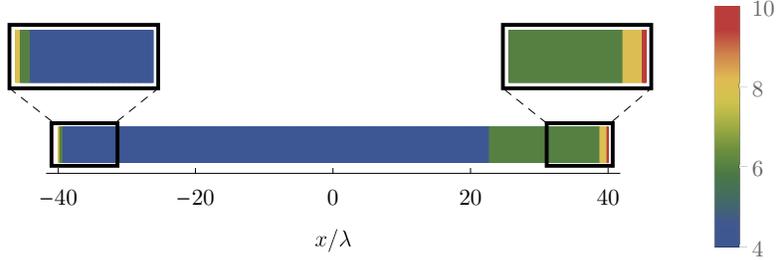}
\end{center}
\caption{%
Order of moment approximation in the computational domain for the goal-adaptive approximation of the solution to the shock-structure problem, including zooms of the regions near the boundaries.
\label{fig:SSvdeg}.}
\end{figure}

\begin{figure}
\centering
\subfloat[Distribution $\beta(g_{\{\IM_{\kappa}\}}^{h,p})$]{\includegraphics[width=0.8\textwidth]{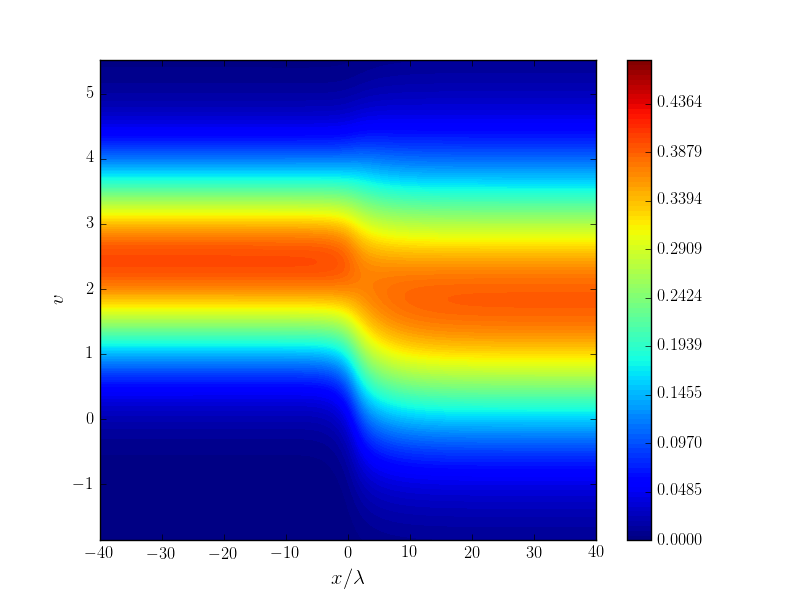}}\,
\subfloat[Dual solution $z^{h,p}_{\{\IM_{\kappa*}\}}$]{\includegraphics[width=0.8\textwidth]{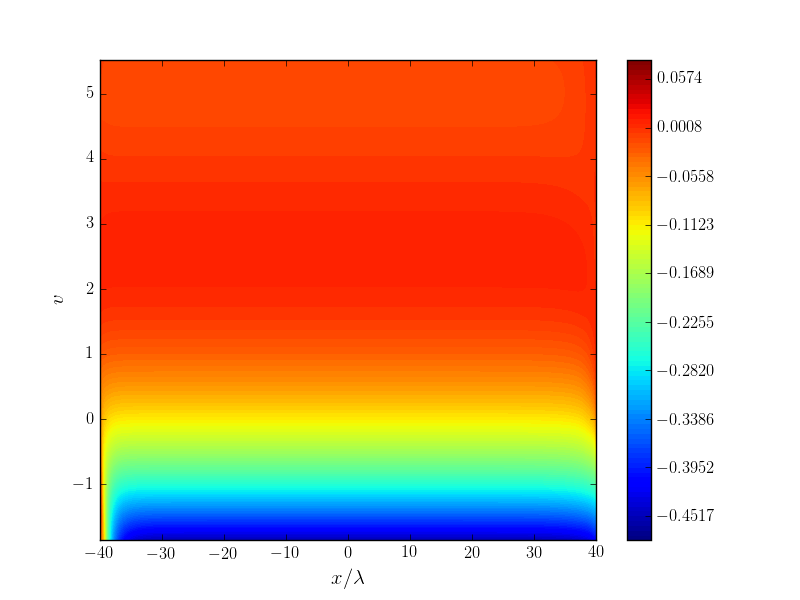}}
\caption{%
Goal-adaptive approximation of the primal distribution $\beta(g_{\{\IM_{\kappa}\}}^{h,p})$ 
for the shock-structure problem ({\it top}) and corresponding approximate dual solution $z^{h,p}_{\{\IM_{\kappa*}\}}$
({\it bottom}) in the final step of the adaptive algorithm\label{fig:SSPD}.}
\end{figure}

The final spatial distribution of the moment orders generated by the goal-adaptive algorithm
is displayed in Figure~\ref{fig:SSvdeg}. One can observe that the goal-adaptive algorithm introduces most of the moment 
refinements near the boundaries and, in particular, near the right boundary. 
To elucidate the refinement pattern, the top and bottom panels in Figure~\ref{fig:SSPD}
display the approximation of the primal solution and of the dual solution, respectively,
in the final step of the adaptive algorithm.
Figure~\ref{fig:SSPD} indicates that the distribution $\beta(g_{\{\IM_{\kappa}\}}^{h,p})$ 
exhibits non-equilibrium behavior in the neighborhood of the shock which is located near the center of the
domain. At further distances from the shock, including the vicinity of the boundary, the distribution 
is close to equilibrium. The dual solution on the other hand manifests boundary layers near the left
and right boundaries.

\begin{figure}
\centering
\includegraphics[width=0.9\textwidth]{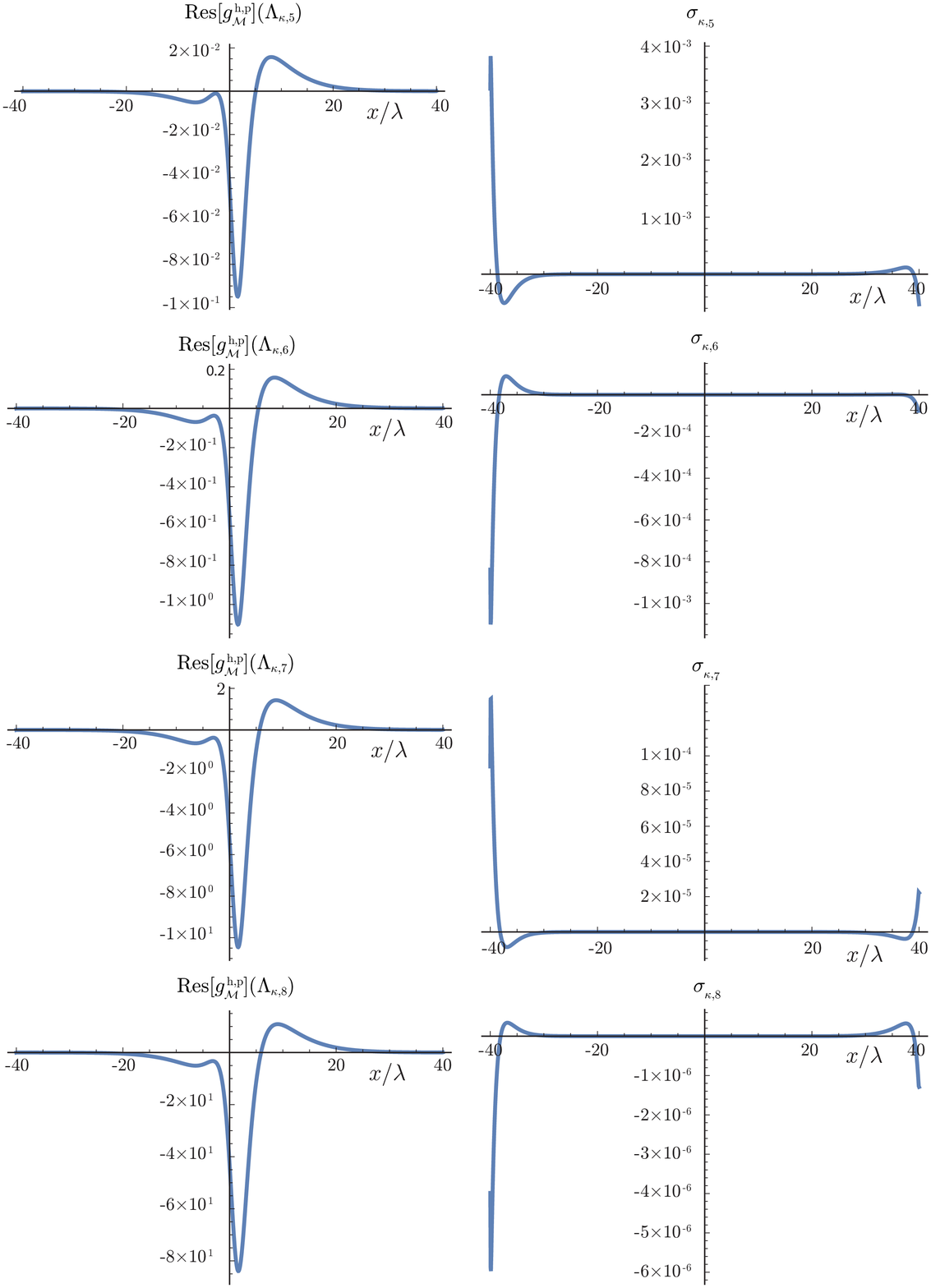}
\caption{Plots of the weigted residuals $\operatorname{Res}[g^{h,p}_{\IM}](\Lambda_{\kappa,i})$ ({\it left\/}) and corresponding dual coefficients $\sigma_{\kappa,i}$ ({\it right\/}) versus the centroids of the 
elements~$\kappa\in\mathcal{T}^h$ in the first step of the goal-adaptive algorithm for the shock-structure 
problem.\label{fig:DWR}}
\end{figure}

\begin{figure}
\centering
\includegraphics[width=0.7\textwidth]{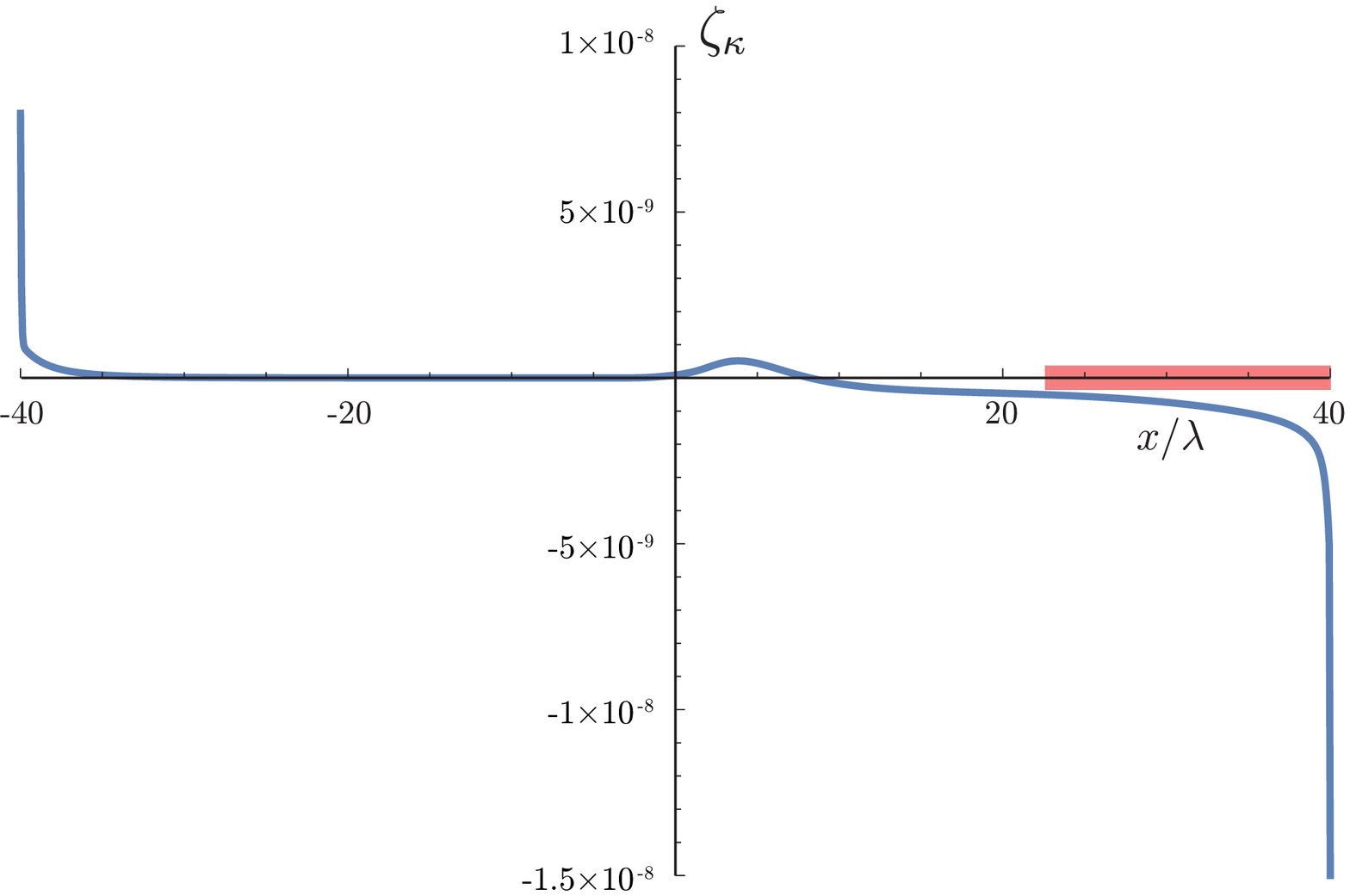}
\caption{Plot of the error indicators $\zeta_{\kappa}$ in~\EQ{zetadef} versus the element centroids in 
the first step of the adaptive algorithm for the shock structure problem. The red interval indicates
the region where moment refinement is introduced.\label{fig:ErrInd}}
\end{figure}

A more detailed view of the origin of the moment-order refinement pattern in Figure~\FIG{SSvdeg} is provided by
Figures~\FIG{DWR} and~\FIG{ErrInd}. Figure~\FIG{DWR} displays the components of the element-wise error
indicators in~\EQ{zetadef}, viz. $\operatorname{Res}[g^{h,p}_{\IM}](\Lambda_{\kappa,i})$ and $\sigma_{\kappa,i}$,
in the initial approximation, i.e. for $M=4$ and $M_*=8$. For the
considered one-dimensional test case and piecewise constant DGFE approximation, the basis functions 
$\Lambda_{\kappa,i}(x,v)$ correspond to monomials in velocity dependence of order~$i$ supported on the
element $\kappa\in\mathcal{T}^h$:
\begin{equation}
\label{eq:Lambda}
\Lambda_{\kappa,i}(x,v)=
\begin{cases}
v^i&\quad\text{if }x\in\kappa
\\
0&\quad\text{otherwise}
\end{cases}
\qquad
i\in\{0,1,\ldots,M_{\kappa*}\}
\end{equation}
Each $\sigma_{\kappa,i}$ represents the corresponding weight of the approximate dual solution.
Let us note that in Figure~\ref{fig:DWR} we have omitted the terms of order $\leq{}4$ 
because $\operatorname{Res}[g^{h,p}_{\IM}](\Lambda_{\kappa,i})$ vanishes for $i\leq{}4$ on account of
Galerkin orthogonality; see Section~\SEC{ErrEst}. The moments of the residual in Figure~\ref{fig:DWR} 
({\em left\/}) are indeed largest in the region where the largest deviations from equilibrium occur, viz. in 
the neighborhood of the shock. The coefficients of the dual solution in Figure~\ref{fig:DWR} 
({\em right\/}) however indicate that the contribution 
of this region to the quantity of interest is negligible. Instead, the goal functional is most sensitive to 
errors in the neighborhood of the boundaries. Multiplication of the weighted residuals 
$\operatorname{Res}[g^{h,p}_{\IM}](\Lambda_{\kappa,i})$ and the dual coefficients $\sigma_{\kappa,i}$
and summation within each element yields the element-wise error indicators $\zeta_{\kappa}$ as depicted
in Figure~\FIG{ErrInd}. Figure~\FIG{ErrInd} indicates that despite the fact that non-equilibrium effects are
most prominent in the center of the domain near the shock, the largest contribution to the error in the
goal quantity originates near the boundaries. The elements in the vicinity of the boundary thus qualify for
refinement. The red interval in Figure~\ref{fig:ErrInd} indicates the region that is marked for refinement
after cancellations have been accounted for.

\section{Conclusion}
\label{sec:Conc}
In this work we introduced a new goal-oriented a-posteriori error analysis and an adaptive-refinement strategy for numerical 
approximation of the steady Boltzmann equation. The approximation is  based on a combination of moment-system approximation in velocity dependence and discontinuous-Galerkin finite-element approximation in spatial dependence. 
We considered a moment-closure relation derived from the minimization of a divergence-based relative entropy. 
The combined DGFE moment method can be construed as a Galerkin finite-element approximation of the Boltzmann equation in renormalized form, based on a tensor-product approximation space composed of the DGFE approximation space in position dependence and global polynomials in velocity dependence. We introduced a numerical flux for the DGFE scheme based on the
position-velocity upwind distribution in the DGFE moment approximation.

The goal-oriented a-posteriori error estimate that we considered is of the usual dual-weighted residual form, furnished with
a linearized dual problem. By virtue of the selected upwind-distribution-based numerical flux, the prerequisite 
linearization is straightforward independent of the moment order. To enhance the efficiency of the adaptive algorithm, 
we introduced a marking strategy that accounts for cancellations of error contributions between elements,
as opposed to the conventional marking strategies based on error bounds derived from the triangle inequality.
The refinement strategy in the adaptive algorithm is based on local, element-wise increments of the moment-system order.
The proposed adaptive strategy for the Boltzmann equation exploits the Galerkin form of the DGFE moment method
and the hierarchical character of the moment-system approximation.

We presented numerical results for two one-dimensional test cases, viz. a heat-transfer and a shock-structure problem. 
For these test cases we considered a goal functional corresponding to the heat flux. We generally observed good 
agreement between the goal-oriented error estimate and the actual error. Moreover, the proposed upper bound that 
accounts for cancellations was found to be sharp relative to the error estimate, in contrast to the standard
triangle-inequality-based bound. The numerical results demonstrate that the goal-adaptive refinement procedure
provides a highly efficient approximation of the quantity of interest, relative to uniform moment refinement. 

The proposed adaptive moment method can be interpreted as a heterogeneous multiscale method (HMM) of 
type A that introduces a domain decomposition into regions where models of different levels of sophistication
are applied, where the various models corresponding to different members of the moment-system hierarchy.  
The goal-oriented adaptive-refinement strategy performs the domain decomposition and the selection of the
local models in a fully automated and optimal manner.

\section*{Acknowledgment} The authors thank Gertjan van Zwieten (Evalf Computing) for his support of the 
software implementations in the Nutils library (\url{www.nutils.org}).

\bibliography{BibFile}   
\end{document}